\documentclass[11pt]{amsart}
\usepackage{latexsym, amsmath, amssymb,amsthm,amsopn,amsfonts}
\usepackage{version}
\usepackage{epsfig,graphics,color,graphicx,graphpap}
\usepackage{amssymb}
\usepackage{todonotes}
\usepackage{listings}
\usepackage{mathrsfs}
\usepackage[citecolor=blue]{hyperref}
\usepackage{soul,xcolor}

\usepackage[framemethod=tikz]{mdframed}
\newcommand{\redsout}{\bgroup\markoverwith{\textcolor{red}{\rule[0.5ex]{2pt}{.4pt}}}\ULon}
\numberwithin{equation}{section}
\newtheorem{theorem}{Theorem}[section]
\newtheorem*{theorem*}{Theorem}

\newtheorem{proposition}[theorem]{Proposition}
\newtheorem*{theorem_0*}{Proposition}
\newtheorem{lemma}[theorem]{Lemma}

\newtheorem{remark}{Remark}[section]

\title[]{Stable Determination of Coefficients in Nonlinear Dynamical Schr\"{o}dinger Equations by Carleman Estimates}
\author[P. Arrepu]{Pranav Arrepu}
\address{Department of Mathematics, University of California Santa Barbara, Santa Barbara, CA 93106, USA}
\email{pranav\_arrepu@ucsb.edu}

\author[H. Zhou]{Hanming Zhou}
\address{Department of Mathematics, University of California Santa Barbara, Santa Barbara, CA 93106, USA}
\email{hzhou@math.ucsb.edu}

\begin{document}

\begin{abstract}
    We consider the inverse problem of recovering stationary coefficients in a class of dynamical Schr\"odinger equations with locally analytic nonlinear terms. Upon treating the well-posedness for small initial data and trivial boundary data, we proceed to establish stable and unique determination provided knowledge of the coefficients near the boundary and the measured Neumann data of the solution. We discuss both the case of measurement on a subset of the boundary satisfying a certain geometric condition and, under stronger assumptions on the regularity and size of the coefficients, the case of measurement on arbitrary subsets of the boundary. Our argument relies on high-order linearization and Carleman estimates for the linear Schr\"odinger equation. %derived by L. Baudouin and J.-P. Puel (Inverse Problems \textbf{18} (2002) 1537-1554).
\end{abstract}

\maketitle
%\tableofcontents

\section{Introduction}\label{Introduction}

Letting $\Omega$ be a smooth bounded domain in $\mathbb R^{n}\,(n = 2, 3)$, we consider the following nonlinear Schr\"{o}dinger equation: 
\begin{equation} \label{quadraticStartEqn} 
\begin{cases}\left(i\partial_t + \Delta + p(x)\right) u(x,t) + q(x) N(u(x,t),\overline{u(x,t)}) = 0 \text{ on } \Omega \times (0, T),\\ u(x,t) = 0 \text{ on }\partial \Omega \times (0,T),\\ u(x,0) = f(x) \text{ on } \Omega \times \{0\},
\end{cases}
\end{equation}
where $N(u,\overline u)$ represents a general nonlinear term satisfying certain conditions that will be specified. Here $p(x)$ and $q(x)$ are the linear and nonlinear coefficients respectively, and $f(x)$ defines the initial condition.\\ \indent
We briefly record some notations used throughout the paper. We adopt the convention that $L^{2}(\Omega)$ and its subspaces are complex Hilbert spaces. In particular, $L^{2}(\Omega)$ is equipped with the Hermitian inner product
\begin{equation}
\left(f, g \right)_{L^{2}(\Omega)} = \int_{\Omega}f\overline{g}\,dx.
\end{equation}
For $k \geq 0$ an integer, we consider the Sobolev space of complex-valued functions $H^{k}(\Omega)$ with norm given by
\begin{equation} \label{Sobolev_norm_definition}
\left \Vert \cdot \right \Vert_{H^{k}(\Omega)} = \left(\sum_{|\alpha| = 0}^{k}\left\Vert D^{\alpha}\left( \cdot \right)\right\Vert^{2}_{L^{2}(\Omega)}\right)^{\frac{1}{2}}
\end{equation}
and $H^{k}_{0}(\Omega)$, the completion of smooth compactly supported test functions $C^{\infty}_{c}(\Omega)$ with respect to the $H^{k}(\Omega)$-norm. Whenever convenient, we denote $D(A) = H^2(\Omega)\cap H^{1}_{0}(\Omega)$ with norm given by $\Vert \cdot \Vert_{D(A)} = \Vert \cdot \Vert_{H^{2}(\Omega)}$ and $C^{0}_{D(A)} = C([0,T];\,D(A))$ the complex Banach space with norm given by
\begin{equation}
\left \Vert \cdot \right \Vert_{C^{0}_{D(A)}} = \sup_{t\in[0,T]} \left \Vert \cdot \right \Vert_{D(A)}.
\end{equation}
We will also use the fact that in $\mathbb{R}^{n}\,(n = 2, 3)$, $H^{2}(\Omega)$ (and consequently $D(A)$) is a Banach algebra. That is, there exists a constant $K^{*} = K^{*}(\Omega) > 0$ such that for any $h_{1}, h_{2} \in H^{2}(\Omega)$,
\begin{equation} \label{Banach_algebra_constant_0}
\Vert h_1h_2\Vert_{H^{2}(\Omega)} \leq K^{*}\Vert h_1\Vert_{H^{2}(\Omega)}\Vert h_2\Vert_{H^{2}(\Omega)}.
\end{equation}
We now delineate our assumptions on the nonlinearity $N(u,\overline{u})$ appearing in \eqref{quadraticStartEqn}. For some integer $k > 1$ and another constant $\delta > 0$, we suppose that $N = N(z_1, z_2):\mathbb{C} \times \mathbb{C} \rightarrow \mathbb{C}$ is a mapping that is analytic on the disk
\begin{equation} \label{nonlinear_term_disk_of_analyticity}
B_{\delta}(0,0) = \left\{(z_1, z_2) \in \mathbb{C} \times \mathbb{C}\,\Big|\, |z_1|^{2} + |z_2|^{2} < \delta \right\}
\end{equation}
and which satisfies the conditions that
\begin{equation} \label{nonlinearity_condition_0}
\sum_{m = 0}^{k}\binom{k}{m}\cdot \partial^{m}_{z_{1}}\partial^{k - m}_{z_{2}}N(0,0) \neq 0,\;
\partial^{m}_{z_{1}}\partial^{k - m}_{z_{2}}N(0,0) \in \mathbb{R} \text{ for } 0 \leq m \leq k.
\end{equation}
\begin{remark}
In the present study, we use Carleman estimates based on the interval $(-T, T)$ and the assumption that the terms above are real-valued ensures that the solution to a linearized version of \eqref{quadraticStartEqn}, namely \eqref{Partial_t_Second_Order_Difference}, takes a simple form on $(-T,T)$. That the above sum is non-zero ensures that the initial condition of \eqref{Partial_t_Second_Order_Difference} is non-zero.
\end{remark}
We further stipulate for $(h_1,h_2) \in C^{0}_{D(A)} \times C^{0}_{D(A)}$ that
\begin{equation} \label{nonlinearity_conditions}
\begin{cases}
\left \Vert N(h_1,h_2)\right \Vert_{C^{0}_{D(A)}} \leq C_{0}\left \Vert h_1 \right \Vert_{C^{0}_{D(A)}}^{m_{0}}\left \Vert h_2 \right \Vert_{C^{0}_{D(A)}}^{n_{0}},\\
\partial^{\alpha_1}_{z_1}\left(\partial^{\alpha_2}_{z_2} N(0,0)\right) = 
\partial^{\alpha_2}_{z_2}\left(\partial^{\alpha_1}_{z_1} N(0,0)\right) = 0 \text{ for }\alpha_1 + \alpha_2 \leq k - 1,
\end{cases}
\end{equation}
where $C_{0} > 0$ and $m_0 + n_0 > 1$. The condition that $m_{0} + n_{0} > 1$ is required for the well-posedness and linearization of \eqref{quadraticStartEqn}. Also, we point out that the condition in \eqref{nonlinearity_conditions} only needs to hold in some neighborhood of $(0,0) \in C^{0}_{D(A)} \times C^{0}_{D(A)}$.\\ \indent
In fact, the analyticity of $N(z_{1},z_{2})$ and \eqref{nonlinearity_conditions} imply that when
\begin{equation} \label{absolute_convergence_disk}
\Vert h_{1} \Vert_{C^{0}_{D(A)}}^{2} + \Vert h_{2} \Vert_{C^{0}_{D(A)}}^{2} < \frac{\delta}{\left(K^{*}\right)^2},
\end{equation}
$N(h_{1}, h_{2})$ is given by the power series
\begin{equation} \label{nonlinear_term_power_series}
\sum_{l = k}^{\infty}\left(\sum_{m = 0}^{l} \frac{1}{m!(l-m)!} \left(\partial_{z_1}^{l-m}\partial_{z_2}^{m}N(0,0) \right)h_{1}^{l-m}h_{2}^{m} \right)
\end{equation}
converging absolutely with respect to $\Vert \cdot \Vert_{C^{0}_{D(A)} \times C^{0}_{D(A)}}$. To see this, begin by observing that as a result of \eqref{Banach_algebra_constant_0}, we have
\begin{equation} \label{nonlinear_term_power_series_1}
\begin{split}
\left \Vert h_{1}^{l-m}h_{2}^{m} \right \Vert_{C^{0}_{D(A)}} \leq \left( K^{*} \right)^{l-1}\cdot \left\Vert h_{1} \right\Vert^{l-m}\left\Vert h_{2} \right\Vert^{m}
\end{split}
\end{equation}
so that
\begin{equation}
\begin{split}
K^{*}\left \Vert h_{1}^{l-m}h_{2}^{m} \right \Vert_{C^{0}_{D(A)}} \leq \left( K^{*} \right)^{l}\cdot \left\Vert h_{1} \right\Vert^{l-m}\left\Vert h_{2} \right\Vert^{m} = \left\Vert K^{*}h_{1} \right\Vert^{l-m}\left\Vert K^{*}h_{2} \right\Vert^{m}.
\end{split}
\end{equation}
Now, \eqref{nonlinear_term_power_series} follows since the analyticity of $N(z_{1},z_{2})$ on the disk $B_{\delta}(0,0)$ implies that for $(h_1,h_2) \in C^{0}_{D(A)} \times C^{0}_{D(A)}$ satisfying \eqref{absolute_convergence_disk},
\begin{equation} \label{nonlinear_term_power_series_2}
\sum_{l = k}^{\infty}\left(\sum_{m = 0}^{l} \frac{1}{m!(l-m)!} \left|\partial_{z_1}^{l-m}\partial_{z_2}^{m}N(0,0) \right|\left \Vert K^{*}h_{1}\right \Vert^{l-m}\left \Vert K^{*}h_{2} \right \Vert^{m} \right) < \infty,
\end{equation}
where the first $k-1$ terms vanish by \eqref{nonlinearity_conditions}.\\ \indent
The nonlinearity in \eqref{quadraticStartEqn} is motivated by terms appearing in the Gross-Pitaevskii equation\footnote{See \cite{GPEquation} for further details and applications to the modelling of Bose-Einstein condensates.}, where $N(u,\overline u)=|u|^2 u$. In this case, the linear coefficient $p$ represents an external potential while the nonlinear term models particle interactions with the coupling constant $q$. Then the solution $u$ defines the wavefunction. We remark that the condition in \eqref{nonlinearity_conditions} is similar to a growth condition satisfied by the nonlinear terms studied for a transport equation in \cite{SemilinearTransportNonlinearRecovery}. We also point out that non-polynomial, analytic nonlinear terms have been considered in the inverse source problem for certain elliptic equations in \cite{SemilinearEllipticSource}. \\ \indent
For \eqref{quadraticStartEqn}, we establish the following well-posedness result when the initial data $f$ is small:
\begin{theorem} \label{nonlinear_well-posedness}
Let $\Omega$ be a smooth bounded domain in $\mathbb{R}^n\; (n = 2, 3)$. Let $p(x) \in L^{\infty}(\Omega)$ be real-valued. Consider $q(x) \in H^{2}(\Omega), f(x) \in D(A)$ and suppose $N(z_1,z_2)$ is an analytic mapping near $(0,0)$ that satisfies \eqref{nonlinearity_condition_0}-\eqref{nonlinearity_conditions}. Then, there exists $r > 0$ sufficiently small so that when $\Vert f \Vert_{D(A)} \leq r$, there is a unique solution $u$ to \eqref{quadraticStartEqn} in $C([0,T];D(A)) \cap C^{1}([0,T];L^2(\Omega))$ satisfying 
$$\Vert u \Vert_{C^{0}_{D(A)}} \leq C\Vert f\Vert_{D(A)}$$
for some constant $C>1$ only depending on $\Omega$ and $\|p\|_{L^\infty(\Omega)}$.
\end{theorem}
\begin{remark}
\hfill
\begin{itemize}
\item[(a)] We apply the Banach algebra property in \eqref{Banach_algebra_constant_0} to prove the above result and this restricts our consideration to $\mathbb{R}^{n}\; (n = 2, 3)$. By assuming higher regularity for $q(x)$ and $f(x)$, we can establish a similar result in $\mathbb{R}^{n}$ for $n > 3$.
\item[(b)] We note that the well-posedness in $C([0,T];D(A)) \cap C^{1}([0,T];L^2(\Omega))$ of equations similar to \eqref{quadraticStartEqn}, namely those of the form ``$i\partial_{t}u + \Delta u = F(u)$" when $F$ is locally Lipschitz, has been treated in the previous literature (see Lemma 4 in \cite{Brezis_Gallouet_NLSE} and Theorem 1 in \cite{Segal_Nonlinear_Semigroups}) by Picard iterations. However, Theorem \ref{nonlinear_well-posedness} is of interest since in certain cases, it provides a sharper estimate of the solution than what would be guaranteed by using the results of \cite{Brezis_Gallouet_NLSE,Segal_Nonlinear_Semigroups} and applying Gronwall's Inequality to obtain a stability estimate. As will be shown in the proof of the above theorem, this is a consequence of passing the dependence on $T > 0$ onto the size of the initial data $f(x)$.
\end{itemize}
\end{remark}
We are interested in the inverse problem of uniquely and stably determining $p(x)$ and $q(x)$ from the Neumann data measurement $(\partial u/\partial \nu)\big \vert_{\partial \Omega \times (0,T)}$. This is inspired by the well-known Calder\'{o}n problem for the conductivity equation ``$\nabla \cdot (\gamma \nabla u) = 0$", which seeks to determine the coefficient $\gamma$ from the measurements $f \mapsto \gamma (\partial u/\partial \nu)\big|_{\partial \Omega}$ given Dirichlet data $u = f$ on $\partial \Omega$. Since the seminal work of \cite{SU1987}, the Calder\'{o}n problem has been studied in numerous variations. We briefly mention references for previous work on inverse problems for the Schr\"odinger equation. When $q\equiv 0$, \eqref{quadraticStartEqn} becomes a linear Schr\"odinger equation. The inverse coefficient problem for linear Schr\"odinger equations has been extensively studied for both stationary coefficients \cite{StartingCarlemanEstimate, MainMethod, BKS2016, Eskin2003, HKSY2019, KSS2014, MOR2008,YY2010} and time-dependent coefficients \cite{BB2020, CKS2015, Eskin2008, GalerkinKianSoccorsi}; see also the references therein. The studies in \cite{BB2020, MainMethod} consider the partial data inverse problem of stable and unique determination from a measurement of the Neumann data on an arbitrary subset of $\partial \Omega$. \\ \indent
Fewer results are known when nonlinear terms appear. When $N(u,\overline u)$ is of the form
\begin{equation} \label{polynomial_special_case}
N(u,\overline{u}) = \sum_{a + b = 2}^{K} q_{a,b} u^{a}\overline{u}^{b},
\end{equation}
uniqueness and stability for the recovery of time-dependent coefficients were established in \cite{PDNSE, MagneticGPEqn} through the construction of geometric optics solutions. When $N(u,\overline u)=u^2$, the work \cite{TimeDependentCase} considers the related inverse problem of determining both linear and nonlinear coefficients from the source-to-solution map. Related to \eqref{polynomial_special_case} is the recovery of coefficients when $N(u,\overline{u})$ is of the form
\begin{equation} \label{polynomial_special_case_0}
N(u,\overline{u}) = \sum_{l = 2}^{\infty}\left(\sum_{a + b = l} q_{a,b}u^{a}\overline{u}^{b}\right).
\end{equation}
Since each term $q_{a,b}$ is the coefficient of a finite-degree polynomial just as in \eqref{polynomial_special_case}, the recovery of the nonlinear term in \eqref{polynomial_special_case_0} can be treated by induction on $l \in \mathbb{N}$ using the methods in \cite{MagneticGPEqn}. Similar induction arguments are used in \cite{KU2020} for recovering coefficients in elliptic equations with a nonlinear term resembling \eqref{polynomial_special_case_0} and likewise in \cite{SemilinearTransportNonlinearRecovery} for nonlinear transport equations. \\ \indent
The main contribution of the present paper is the use of linearization to recover the coefficient of a term $N(u, \overline{u})$ satisfying \eqref{nonlinearity_condition_0}-\eqref{nonlinearity_conditions}, which generalizes the nonlinearities considered in previous works by permitting non-polynomial terms. We introduce a small parameter $\epsilon \in \mathbb{R}$ in the initial condition of \eqref{quadraticStartEqn} and consider
\begin{equation} \label{quadraticStartEqn_0} 
\begin{cases}\left(i\partial_t  + \Delta  + p_{j}(x) \right)u_{\epsilon,j}(x,t) + q_{j}(x) N(u_{\epsilon,j}(x,t),\overline{u_{\epsilon,j}(x,t)}) = 0 \text{ on } \Omega \times (0, T),\\ u_{\epsilon,j}(x,t) = 0 \text{ on }\partial \Omega \times (0,T),\\ u_{\epsilon,j}(x,0) = \epsilon f(x) \text{ on } \Omega \times \{0\},
\end{cases}
\end{equation}
for $j = 1,2$. Note that the well-posedness of \eqref{quadraticStartEqn_0} for $\epsilon \in \mathbb{R}$ small follows from Theorem \ref{nonlinear_well-posedness}.\\ \indent
From a physical point of view, our inverse problem consists in determining both the linear and nonlinear potentials of a medium by probing it with small disturbances generated at the initial time. Our data are the responses of the medium to these small disturbances measured on part of the boundary.

For a subset $S \subseteq \partial \Omega$ and $l \in \mathbb{N}$, define
\begin{equation}
\delta^{(l)}_{S} = \left\Vert \partial_t\left(\partial^{l}_{\epsilon}\partial_{\nu}(u_{\epsilon,1}-u_{\epsilon,2})\Big|_{\epsilon = 0}\right) \right\Vert_{L^{2}(S \times (0, T))}.
\end{equation}
We list our main results as follows. The first result concerns measurement of the Neumann data on a subset of the boundary satisfying a geometric condition.
\begin{theorem} \label{main_theorem}
Consider \eqref{quadraticStartEqn_0} and let $\Omega$ be a bounded domain in $\mathbb{R}^{n}\; (n = 2, 3)$ with smooth boundary. Additionally, let $\omega \subset \Omega$ be a relatively open neighborhood of $\partial \Omega$. 
We take real-valued $f(x) \in H^4(\Omega) \cap H^{3}_{0}(\Omega)$, real-valued $p_j(x) \in H^{2}(\Omega) \; (j = 1, 2)$, and $q_{j}(x) \in H^{2}(\Omega)\; (j = 1, 2)$ such that $q_{1}(x) - q_{2}(x)$ is real-valued. 
Suppose that there exist constants $M, \gamma_{-}, \gamma_{+} > 0$ satisfying
\begin{equation} \label{large_boundary_coefficient_constraints}
\begin{cases}
\Vert p_j \Vert_{H^{2}(\Omega)} \leq M\; (j = 1, 2),\\
|f(x)| \geq \gamma_{-} \text{ on } \Omega \setminus \omega,\\
\Vert f \Vert_{H^{4}(\Omega)} \leq \gamma_{+}.
\end{cases}
\end{equation}
Further, let $N(z_1,z_2)$ be an analytic mapping near $(0,0)$ that satisfies \eqref{nonlinearity_condition_0}-\eqref{nonlinearity_conditions} and for some $x_{0} \notin \overline{\Omega}$, consider a subset $\Gamma_{0} \subseteq \partial \Omega$ such that
\begin{equation} \label{large_boundary_definition}
\left\{x \in \partial\Omega\,\big|\, (x - x_{0})\cdot \nu(x) \geq 0 \right\} \subseteq \Gamma_{0}.
\end{equation}
\begin{itemize}
\item[(a)] Suppose that $p_1(x) = p_2(x)$ on $\omega \subset \Omega$. Then, one has that
\begin{equation} \label{main_theorem_stable_determination_0}
    \Vert p_{1} - p_{2}\Vert_{L^2(\Omega \setminus \omega)} \leq C \delta^{(1)}_{\Gamma_{0}}
\end{equation}
for some constant $C = C(\Omega, \Gamma_{0}, M, T, \gamma_{-}, \gamma_{+}) > 0$. 
\item[(b)] Suppose instead that $p_1(x) = p_2(x)$ on $\Omega$ and that $q_1(x) = q_2(x)$ on $\omega \subset \Omega$. Then, one has that
\begin{equation} \label{main_theorem_stable_determination}
    \Vert q_{1} - q_{2}\Vert_{L^2(\Omega \setminus \omega)} \leq C \delta^{(k)}_{\Gamma_{0}}
\end{equation}
for some constant $C = C(\Omega, \Gamma_{0}, M, T, \gamma_{-}, \gamma_{+}) > 0$. 
\item[(c)] Lastly, suppose that $p_1(x) = p_2(x)$ and $q_1(x) = q_2(x)$ on $\omega \subset \Omega$. If there exists $\epsilon_{0} > 0$ such that $(\partial u_{\epsilon,1}/\partial \nu)\vert_{\Gamma_{0}} = (\partial u_{\epsilon,2}/\partial \nu)\vert_{\Gamma_{0}}$ for all $\epsilon \in (0, \epsilon_{0})$, then $p_1(x) = p_2(x)$ and $q_1(x) = q_2(x)$ a.e. on $\Omega \setminus \omega$.
\end{itemize}
\end{theorem}
In the case of measurement on an arbitrary relatively open subset of the boundary, we have the following result.
\begin{theorem} \label{main_theorem_partial_data}
Consider \eqref{quadraticStartEqn_0} and let $\Omega$ be a bounded domain in $\mathbb{R}^{n}\; (n = 2, 3)$ with smooth boundary. Additionally, let $\Gamma$ be a relatively open subset of $\partial \Omega$ and $\omega \subset \Omega$ a relatively open neighborhood of $\partial \Omega$ such that $\partial\omega \setminus \partial\Omega$ is $C^{2}$.\\ \indent
Let $f(x)$, $p_{j}(x)$, and $q_{j}(x)\; (j = 1,2)$ satisfy the conditions in Theorem \ref{main_theorem}. Assume also that $p_{1}(x) - p_{2}(x),\, q_{1}(x) - q_{2}(x) \in H^{4}(\Omega)$ and that there exists a constant $\tilde{M} > 0$ such that $\Vert p_{1} - p_{2}\Vert_{H^{4}(\Omega)}, \Vert q_{1} - q_{2}\Vert_{H^{4}(\Omega)} \leq \tilde{M}$.
\begin{itemize}
\item[(a)] Suppose that $p_1(x) = p_2(x)$ on $\omega \subset \Omega$. Then, one has that
\begin{equation} \label{main_theorem_partial_stable_determination_0}
    \Vert p_1 - p_2\Vert_{L^2(\Omega \setminus \omega)} \leq C \left[\left|\ln \delta^{(1)}_{\Gamma}\right|^{-1}  +  \delta^{(1)}_{\Gamma} \right]^{1/2}
\end{equation}
for some constant $C = C(\Omega, \Gamma, \omega, M, \tilde{M}, T, \gamma_{-}, \gamma_{+}) > 0$. 
\item[(b)] Now assume that $p_1(x) = p_2(x)$ on $\Omega$ and suppose $q_1(x) = q_2(x)$ on $\omega \subset \Omega$. Then, one has that
\begin{equation} \label{main_theorem_partial_stable_determination}
\begin{split}
    \Vert q_1 - q_2\Vert_{L^2(\Omega \setminus \omega)} \leq C\left[\left|\ln \delta^{(k)}_{\Gamma}\right|^{-1}  +  \delta^{(k)}_{\Gamma} \right]^{1/2}
\end{split}
\end{equation}
for some constant $C = C(\Omega, \Gamma, \omega, M, \tilde{M}, T, \gamma_{-}, \gamma_{+}) > 0$.
\item[(c)] Lastly, suppose that $p_1(x) = p_2(x)$ and $q_1(x) = q_2(x)$ on $\omega \subset \Omega$. If there exists $\epsilon_{0} > 0$ such that $(\partial u_{\epsilon,1}/\partial \nu)\vert_{\Gamma_{0}} = (\partial u_{\epsilon,2}/\partial \nu)\vert_{\Gamma_{0}}$ for all $\epsilon \in (0, \epsilon_{0})$, then $p_1(x) = p_2(x)$ and $q_1(x) = q_2(x)$ a.e. on $\Omega \setminus \omega$.
\end{itemize}
\end{theorem}
\begin{remark} \label{inverse_problem_setting_Sobolev_embedding}
\hfill
\begin{itemize}
\item[(a)] The assumption that $\partial \omega \setminus \partial \Omega$ is $C^{2}$ is needed for the proof of Proposition \ref{parabolic_weight_function_0}, which introduces a function used in a parabolic Carleman estimate involved in the proof of the above theorem. We refer to Remark \ref{smoothness_Pi_inner_boundary} for further details. However, this regularity assumption is not restrictive since if the boundary is rough, we can replace $\omega$ by a smaller neighborhood $\omega'$ such that $\partial \omega' \setminus \partial \Omega$ is $C^{2}$ and the lower bound $\gamma_{-} > 0$ in \eqref{large_boundary_coefficient_constraints} with a possibly smaller constant $\gamma'_{-} > 0$.
\item[(b)] Aside from ensuring well-posedness by Theorem \ref{nonlinear_well-posedness}, we consider the inverse problem in $\mathbb{R}^{n}\; (n = 2, 3)$ in order to apply the Sobolev Embedding Theorem in the proofs of Theorem \ref{main_theorem}-\ref{main_theorem_partial_data}. By assuming higher regularity for the coefficients $p_{j}(x),\, q_{j}(x)\;(j = 1,2)$ and initial condition $f(x)$ appearing in \eqref{quadraticStartEqn_0}, one can establish similar results in $\mathbb{R}^{n}$ for $n > 3$.
\item[(c)] The previous studies \cite{PDNSE,MagneticGPEqn} use assumptions similar to our requirement that $p_{1}(x) = p_{2}(x)$ and $q_{1}(x) = q_{2}(x)$ for $x \in \omega$ in order to apply a unique continuation estimate from Lemma $2.2$ in \cite{BB2020}. In the case of the non-dynamical linear Schr\"odinger equation with trivial Dirichlet data outside a relatively open subset $\Gamma \subseteq \partial \Omega$, it is known (cf. Exercise $5.7.4$ in \cite{Is2017}) that knowledge of the linear coefficient in a neighborhood of $\partial \Omega$ along with the Dirichlet-to-Neumann map on $\Gamma$ determines the Dirichlet-to-Neumann map on the entire boundary. However, we are concerned with proving logarithmic stability with respect to Neumann data measured only on $\Gamma$.
\end{itemize}
\end{remark}
The present paper deals with the nonlinearity $N(u,\overline{u})$ through Carleman estimates in the style of A. L. Bukhgeim and M. V. Klibanov in \cite{BK1981}. Carleman estimates for the dynamical Schr\"odinger equation were derived for the first time in \cite{StartingCarlemanEstimate} when the Neumann data was measured on a subset of the boundary satisfying a geometric condition. For earlier literature concerning the application of Carleman estimates to inverse problems, we refer the reader to \cite{K1992} and the references therein. Studies subsequent to \cite{StartingCarlemanEstimate} have applied Carleman estimates to inverse problems for the linear Schr\"odinger equation in more general settings. The authors of \cite{MOR2008} generalized the Carleman estimates in \cite{StartingCarlemanEstimate} to enable stability by observation over arbitrary subsets of $\partial \Omega$ provided the existence of a Carleman weight satisfying a weak pseudoconvexity condition on the domain. The pseudoconvexity requirement was relaxed in \cite{MainMethod} by transforming the Schr\"odinger equation into a parabolic equation for which partial data Carleman estimates can be derived. Carleman estimates in negative order Sobolev spaces are derived in \cite{YY2010} and are used to study inverse problems for $L^{p}$-coefficients for $p > 1$. A Carleman estimate for the linear Schr\"odinger equation on an infinite cylindrical domain was derived in \cite{KSS2014}. For the magnetic Schr\"odinger equation, we mention that Carleman estimates were derived in \cite{CS2012} to study the recovery of the magnetic field and in \cite{HKSY2019} to recover both the magnetic field and the linear coefficient. \\ \indent
We prove Theorem \ref{main_theorem} by first applying the high-order linearization technique detailed in \cite{FO2020, KU2020, KLU2018, LinearizationTechnique, SU1997} to \eqref{quadraticStartEqn_0}. With the linearized equation, we apply the strategy used in \cite{StartingCarlemanEstimate} to study the inverse problem for the linear Schr\"odinger equation. In particular, the authors of \cite{StartingCarlemanEstimate} prepare an equation in which the quantity we wish to estimate is part of the initial data. Then, an estimate involving the initial data is derived through integration by parts. Lastly, terms in this estimate are bounded above by a Carleman estimate for the Schr\"odinger operator $i\partial_t + \Delta + p(x)$ with $p(x) \in L^{\infty}(\Omega)$. \\ \indent
The proof of Theorem \ref{main_theorem_partial_data} is similar in that it relies on using an initial data estimate and estimating terms through another Carleman estimate for the Schr\"odinger operator. As in \cite{MainMethod}, we further need to estimate the solution to a linearized version of \eqref{quadraticStartEqn_0} on a subset of $\omega$ in terms of the Neumann data on $\Gamma$. This requires a unique continuation estimate. \\ \indent
We note that the use of an estimate involving the initial data corresponds to the setting of \eqref{quadraticStartEqn_0} where we impose non-trivial initial data but trivial Dirichlet data on the boundary. This is also the setting considered in other works based on the Carleman estimate such as \cite{StartingCarlemanEstimate, MainMethod, HKSY2019}. This is in contrast to \cite{PDNSE, MagneticGPEqn}, where non-trivial Dirichlet data but trivial initial data are imposed.\\ \indent
In Section \ref{forward_problem_well-posedness}, we provide the proof of Theorem \ref{nonlinear_well-posedness}. In Section \ref{main_theorem_proof_tools}, we address the differentiability of the solution to \eqref{quadraticStartEqn_0} with respect to $\epsilon \in \mathbb{R}$ and also record the Carleman estimates to be used in the proofs of Theorem \ref{main_theorem} and Theorem \ref{main_theorem_partial_data}, which are discussed in Section \ref{main_theorem_proof}.
\bigskip

{\bf Acknowledgments:} We thank the anonymous referees for their feedback which improved the paper. The authors would like to thank Prof. Ru-Yu Lai for bringing their attention to the problem and for helpful discussions. P. Arrepu additionally thanks Prof. Ru-Yu Lai for hosting a visit to the University of Minnesota, Twin Cities and for discussions concerning inverse problems for the Schr\"odinger equation. P. Arrepu also thanks Prof. Loc Nguyen (University of North Carolina, Charlotte) for pointing out an issue with the parameter dependence of the constants appearing in the statements of the main theorems. The authors are partly supported by the NSF grants DMS-2109116 and DMS-2408369, and Simons Foundation Travel Support for Mathematicians MPS-TSM-00008046.

\section{Well-Posedness for the Forward Problem} \label{forward_problem_well-posedness}
We proceed to establish the well-posedness of \eqref{quadraticStartEqn}. Our argument relies on the well-posedness of the linear Schr\"odinger equation taken up in the next section.
\subsection{The Homogeneous Linear Equation} \label{the_homogeneous_linear_equation}
For $p(x) \in L^{\infty}(\Omega)$, we consider
\begin{equation} \label{main_linear_equation}
\begin{cases} \left(i\partial_t + \Delta + p(x)\right) u(x,t) = 0 \text{ on } \Omega \times (0,T),\\ u(x,t) = 0 \text{ on } \partial \Omega \times (0,T),\\ u(x,0) = f(x) \text{ on } \Omega \times \{0\}. \end{cases}
\end{equation}
In \cite{GalerkinKianSoccorsi}, the authors consider the well-posedness of an equation similar to \eqref{main_linear_equation} and it is possible to use their method to show well-posedness in $C\left([0,T]; D(A)\right) \cap C^{1}\left([0,T];L^2(\Omega)\right)$. However, when there is a non-trivial source term in \eqref{main_linear_equation}, the stability estimate in \cite{GalerkinKianSoccorsi} is such that it estimates the solution in $ C\left([0,T]; D(A)\right) \cap C^{1}([0,T];L^2(\Omega))$ in terms of the source in $H^{1}(0,T; L^{2}(\Omega))$. For the fixed point argument concerning the well-posedness of \eqref{quadraticStartEqn} in Section \ref{the_nonlinear_equation}, it is convenient to have an estimate in which the solution can be estimated by the source in the same norm.\\ \indent
For this purpose, we address the well-posedness of \eqref{main_linear_equation} by noting that the unbounded operator $A = i\Delta + i p(x)$ generates a strongly continuous unitary group on $D(A) = H^2(\Omega) \cap H^{1}_{0}(\Omega) \subset L^{2}(\Omega)$, thereby yielding a solution to \eqref{main_linear_equation} in $C([0,T]; D(A)) \cap C^{1}([0,T];L^2(\Omega))$. A consequence of the following theorem is that when \eqref{main_linear_equation} has a non-trivial source term, we can use Duhamel's Formula (recalled in Section \ref{the_nonlinear_equation}) to estimate the solution by the source term in the same norm. %Such an estimate is well suited for the fixed point argument used in Section \ref{the_nonlinear_equation}.
\begin{theorem} \label{linear_well-posedness}
Let $\Omega$ be a bounded domain in $\mathbb{R}^n$ with smooth boundary $\partial \Omega$. Suppose $p(x) \in L^{\infty}(\Omega)$ is real-valued and $f(x) \in D(A)$ where $D(A) = H^2(\Omega) \cap H^1_{0}(\Omega)$. For a given $T > 0$, there exists a unique solution to \eqref{main_linear_equation} in $C([0,T];D(A)) \cap C^{1}([0,T];L^2(\Omega))$. Moreover, we have for some constant $C_{1} = C_{1}\left(\Omega,\Vert p\Vert_{L^{\infty}(\Omega)} \right) > 1$ the stability estimate
\begin{equation} \label{H2_linear_stability_estimate}
\Vert u(\cdot\,,t) \Vert_{H^{2}(\Omega)} \leq C_{1} \Vert f\Vert_{H^{2}(\Omega)},
\end{equation}
which holds for each $t \in [0, T]$.
\end{theorem}
\begin{proof}
We begin by noting that the unbounded operator $iA = -\Delta - p(x)$ with $D(iA) = D(A) = H^2(\Omega) \cap H^1_{0}(\Omega) \subset L^2(\Omega)$ is self-adjoint. Indeed, the operator is densely defined and symmetric. One can use elliptic regularity estimates to show that it is closed. Now, the fact that $iA$ is self-adjoint follows\footnote{See, for example, the discussion in Section VIII.2 in \cite{Simon_Reed_Functional}.} from the fact that $(iA)^{*}$ is symmetric. To see this, note that since $p(x) \in L^{\infty}(\Omega)$, multiplication by $p(x)$ is a bounded linear operator on $L^{2}(\Omega)$. From Lemma 6.3.2 in \cite{Buhler_Salamon_Functional}, we obtain 
\begin{equation}
\left(iA\right)^{*} = \left(-\Delta\right)^{*} - p(x) = -\Delta - p(x).
\end{equation}
The above uses the assumption that $p(x)$ is real-valued as well as the self-adjointness of $\left(D(A),-\Delta\right)$, which is treated by Proposition 2.6.1 in \cite{SemilinearSchrodinger_Cazenave}. %We also refer the reader to Section 6.3 in \cite{Buhler_Salamon_Functional} and Section 2.1 in \cite{SemilinearEvolutionEquations_Cazenave} for the self-adjointness of $\left(D(A),-\Delta\right)$. 
The symmetry of $\left(iA\right)^{*}$ follows from the density of the test functions in $D(A)$ with respect to the $H^{1}(\Omega)$-norm and the assumption that $p(x)$ is real-valued.\\ \indent
The existence of a strongly continuous unitary group $\{e^{tA}\}_{t \in \mathbb{R}}$ for the evolution equation 
\begin{equation} \label{evolution_equation}
\frac{du(x,\cdot)}{dt} = Au(x,\cdot)
\end{equation}
and thus a unique solution to \eqref{the_homogeneous_linear_equation} in $C([0,T];D(A))$ now follow from Stone's Theorem (e.g. see Theorem 3.2.1 in \cite{Sinha_Sachi_Semigroup}). This solution is also in $C^1([0,T]; L^2(\Omega))$ since we have $\partial_t u(x,t) = i\Delta u(x,t) + i p(x) u(x,t)$ from the equation itself and $i \Delta u(x,t) + i p(x) u(x,t) \in C([0,T]; L^2(\Omega))$.\\ \indent 
%Details for the stability estimate
It remains to show \eqref{H2_linear_stability_estimate}. The fact that $\{e^{tA}\}_{t \in \mathbb{R}}$ is unitary means 
\begin{equation} \label{L2_linear_stability_estimate}
\Vert u(\cdot\,, t)\Vert_{L^{2}(\Omega)} = \Vert e^{tA}f\Vert_{L^{2}(\Omega)} = \Vert f\Vert_{L^{2}(\Omega)}
\end{equation}
for each fixed $t \in [0,T]$ and $f(x) \in D(A)$. Using the fact that the generator $A = i\Delta + ip(x)$ commutes with the operators $\{e^{tA}\}_{t \in \mathbb{R}}$ and \eqref{L2_linear_stability_estimate}, we have
\begin{equation} \label{H2_linear_stability_estimate_step_1}
\begin{split}
\Vert Au(\cdot\,, t)\Vert_{L^{2}(\Omega)} = \Vert Ae^{tA}f\Vert_{L^{2}(\Omega)}
= \Vert e^{tA}Af\Vert_{L^{2}(\Omega)} = \Vert Af\Vert_{L^{2}(\Omega)}\\
\leq \Vert \Delta f\Vert_{L^{2}(\Omega)} + \Vert p\Vert_{L^{\infty}(\Omega)}\Vert f\Vert_{L^{2}(\Omega)}.
\end{split}
\end{equation}
Moreover, upon applying an elliptic regularity estimate, we obtain for some constant $C = C(\Omega) > 0$ the lower bound
\begin{equation} \label{H2_linear_stability_estimate_step_2}
\Vert Au(\cdot\,, t) \Vert_{L^{2}(\Omega)} \geq \frac{1}{C}\left(\Vert u(\cdot\,, t)\Vert_{H^{2}(\Omega)} - \Vert u(\cdot\,, t)\Vert_{L^{2}(\Omega)} \right) = \frac{1}{C}\left(\Vert u(\cdot\,, t)\Vert_{H^{2}(\Omega)} - \Vert f\Vert_{L^{2}(\Omega)} \right),
\end{equation}
where \eqref{L2_linear_stability_estimate} was used for the last equality above. Then, \eqref{H2_linear_stability_estimate_step_1} becomes
\begin{equation} \label{H2_linear_stability_estimate_step_3}
\Vert u(\cdot\,, t)\Vert_{H^{2}(\Omega)} \leq C\left(1 + \Vert p\Vert_{L^{\infty}(\Omega)}\right)\Vert f\Vert_{H^{2}(\Omega)} + \Vert f\Vert_{H^{2}(\Omega)},
\end{equation}
which establishes \eqref{H2_linear_stability_estimate}.
\end{proof}
\begin{remark} \label{complex_real_field_technical_comment}
One could use the self-adjointness of $\left(D(A),-\Delta\right)$ and Duhamel's Formula to deal with the $p(x)$-term in \eqref{main_linear_equation} and establish well-posedness. Moreover, Gronwall's Inequality could be invoked to obtain a stability estimate. This is the approach taken in Chapter $4$ of \cite{SemilinearEvolutionEquations_Cazenave}. However, the resulting stability estimate would depend on $T > 0$ whereas the estimate proved in the preceding theorem is uniform over $[0,T]$.
\end{remark}

\subsection{The Nonlinear Equation} \label{the_nonlinear_equation}
In this section, we prove Theorem \ref{nonlinear_well-posedness}. We shall use Duhamel's Formula for the following IBVP, which is the same as \eqref{main_linear_equation} but with a non-trivial source term:
\begin{equation} \label{main_inhomogeneous_linear_equation}
\begin{cases} \left(i\partial_t + \Delta + p(x) \right)u(x,t) = g(x,t) \text{ on } \Omega \times (0,T),\\ u(x,t) = 0 \text{ on } \partial \Omega \times (0,T),\\ u(x,0) = f(x) \text{ on } \Omega \times \{0\}.\end{cases}
\end{equation}
If \eqref{main_inhomogeneous_linear_equation} has a solution in $C([0,T];D(A)) \cap C^1([0,T];L^2(\Omega))$, it is uniquely defined by Duhamel's Formula as follows:
\begin{equation}\label{var_param_formula}
u(x,t) = e^{tA}f(x) - i \int_{0}^{t} e^{(t-s)A}g(x,s)\,ds.
\end{equation}
Additionally, if $g(x,t) \in W^{1,1}(0,T;L^2(\Omega))$ or $g(x,t) \in L^1(0,T;D(A))$, then $u(x,t)$ defined by \eqref{var_param_formula} belongs to $C([0,T];D(A)) \cap C^1([0,T];L^2(\Omega))$ (e.g. see Section 1.6 of \cite{SemilinearSchrodinger_Cazenave}).\\ \indent
In the following, recall that $C^{0}_{D(A)} = C([0,T];D(A))$ is equipped with the norm given by
\begin{equation} \label{norm_for_contraction_space}
\Vert w\Vert_{C^{0}_{D(A)}} = \sup_{t \in [0,T]}\Vert w(\cdot\,,t) \Vert_{D(A)} = \sup_{t \in [0,T]}\Vert w(\cdot\,,t) \Vert_{H^2(\Omega)}.
\end{equation}
We denote the closed ball in $C^{0}_{D(A)}$ of radius $r > 0$ by
\begin{equation}
\overline{B_{r}(0)} = \left\{w \in C^{0}_{D(A)}\,\Big|\,\Vert w \Vert_{C^{0}_{D(A)}} \leq r\right\}.
\end{equation}
Recall also from Theorem \ref{linear_well-posedness} the unitary group $\left\{e^{tA}\right\}_{t \in \mathbb{R}}$ generated by $A = i\Delta + ip(x)$ and that given $h(x) \in D(A)$, we have for each $t\in [0, T]$ and some constant $C_{1} > 1$ the estimate 
\begin{equation} \label{linear_Schrodinger_unitary_group_estimate}
\Vert e^{tA}h\Vert_{D(A)} \leq C_{1}\Vert h\Vert_{D(A)}.
\end{equation}
\\
\textit{Step 1 (Obtaining a Fixed Point Equation).}\\
Let us begin by obtaining through a standard argument a fixed point equation that will enable one to apply the Contraction Mapping Principle. %(e.g. see Theorem 5.5 in \cite{ODEs}).  We refer the reader to Proposition 2.2 of \cite{PDNSE} and Proposition 3 of \cite{TimeDependentCase} for the general method and notation used in our proof. 
We shall find a solution of the form $u = v + w$, which are defined in \eqref{homog_part_nonlinear_soln} and \eqref{source_part_nonlin_soln}.\\ \indent
First, we consider
\begin{equation} \label{homog_part_nonlinear_soln}
\begin{cases} i\partial_t v(x,t) + \Delta v(x,t) + p(x) v(x,t) = 0 \text{ on } \Omega \times (0,T),\\ v(x,t) = 0 \text{ on } \partial \Omega \times (0,T), \\ v(x,0) = f(x) \text{ on } \Omega \times \{0\}.\end{cases}
\end{equation}
Since $f(x) \in D(A)$, Theorem \ref{linear_well-posedness} shows that the solution to \eqref{homog_part_nonlinear_soln} is unique in $C([0,T];D(A)) \cap C^1([0,T];L^2(\Omega))$.\\ \indent
Now, denote by $S(g(x,t))$ the solution to the following equation with source term $g(x,t)$ and trivial data:
\begin{equation} \label{source_part_nonlin_soln}
\begin{cases} i\partial_t w(x,t) + \Delta w(x,t) + p(x) w(x,t) = g(x,t) \text{ on } \Omega \times (0,T),\\ w(x,t) = 0 \text{ on } \partial \Omega \times (0,T), \\ w(x,0) = 0 \text{ on } \Omega \times \{0\}.\end{cases}
\end{equation}
If $g(x,t) \in C^{0}_{D(A)}$, the solution defined by \eqref{var_param_formula} is unique in $C([0,T];D(A)) \cap C^1([0,T];L^2(\Omega))$.\\ \indent
Now, we see that solving \eqref{quadraticStartEqn} with initial data $f(x) \in D(A)$ amounts to finding a fixed point $w = w(x,t)$ satisfying
\begin{equation} \label{Fixed_point_nonlinear_eqn}
S(-q(x)N(u,\overline{u})) = S(-q(x)N(v + w,\overline{v + w})) = w.
\end{equation}
Recall that we assume $q(x) \in H^{2}(\Omega)$. So if we take $v, w \in C^{0}_{D(A)}$, the Banach algebra property in \eqref{Banach_algebra_constant_0} and the first condition in \eqref{nonlinearity_conditions} imply that 
\begin{equation}
-q(x)N(v + w,\overline{v + w}) \in C^{0}_{D(A)}.
\end{equation}
Consequently, \eqref{var_param_formula} implies
\begin{equation}
S(-q(x)N(v + w,\overline{v + w})) \in C^{0}_{D(A)}.
\end{equation}\\
\textit{Step 2 (Closure of the Solution Map).}\\
 Let us now consider $\Vert f \Vert_{D(A)} \leq r$, which implies 
\begin{equation}
\Vert v\Vert_{C^{0}_{D(A)}} = \Vert e^{tA}f\Vert_{C^{0}_{D(A)}} \leq C_{1}r.
\end{equation}
In order to apply the Contraction Mapping Principle, we will show that when $\Vert f \Vert_{D(A)} \leq r$ and $w \in \overline{B_{r}(0)}$ for $r > 0$ sufficiently small, we have that
\begin{equation} \label{mapping_condition_0}
S\left(-q(x)N(v + w,\overline{v + w})\right) \in \overline{B_{r}(0)}. 
\end{equation}
If we directly apply the solution formula \eqref{var_param_formula} and invoke the unitary group estimate \eqref{linear_Schrodinger_unitary_group_estimate}, we obtain
\begin{align} \label{check_mapping_condition}
\begin{split}
\Vert S\left(-q(x)N(v + w,\overline{v + w})\right) \Vert_{C^{0}_{D(A)}} \\ =
\left\Vert \int_{0}^{t} e^{(t-s)A}\left[-qN(v + w,\overline{v + w})(\cdot\,,s)\right] \,ds \right\Vert_{C^{0}_{D(A)}} \\ \leq 
\sup_{t \in [0,T]}\left(t\sup_{s \in [0,t]} \left \Vert e^{(t-s)A}\left[-qN(v + w,\overline{v + w})(\cdot\,,s)\right] \right \Vert_{D(A)} \right) 
 \\ \leq 
C_{1} TK^{*}\Vert q \Vert_{H^{2}(\Omega)}\left \Vert N(v + w,\overline{v + w}) \right \Vert_{C^{0}_{D(A)}} \\
\leq C_{0} C_{1} T K^{*} \Vert q\Vert_{H^{2}(\Omega)} \left(\Vert v \Vert_{C^{0}_{D(A)}} + \Vert w \Vert_{C^{0}_{D(A)}}\right)^{m_0 + n_0}.
\end{split}
\end{align}
In \eqref{check_mapping_condition}, we use the fact that we have the estimate
\begin{align} \label{Banach_algebra_property}
\left \Vert -qN(v + w,\overline{v + w})(\cdot\,,t)\right \Vert_{D(A)} \leq 
K^{*}\left \Vert q\right\Vert_{H^{2}(\Omega)} \left\Vert N(v + w,\overline{v + w})(\cdot\,,t) \right \Vert_{D(A)}
\end{align}
for each $t \in [0, T]$, which follows from the Banach algebra property noted in \eqref{Banach_algebra_constant_0}. The estimate
\begin{align}
\left \Vert N(v + w,\overline{v + w}) \right \Vert_{C^{0}_{D(A)}} \leq
C_{0} \left(\Vert v \Vert_{C^{0}_{D(A)}} + \Vert w \Vert_{C^{0}_{D(A)}}\right)^{m_0 + n_0}
\end{align}
follows from \eqref{nonlinearity_conditions}. Now, consider $0 < r \leq r_{m_{0},n_{0}}$ where
\begin{equation} \label{stability_constants_nonlinear_well-posedness_0}
r_{m_0,n_0} = \left(\frac{1}{C_{0} C_{1} \left(1 + C_{1}\right)^{m_0 + n_0} T K^{*} \Vert q\Vert_{H^{2}(\Omega)} }\right)^{1/\left(m_{0} + n_{0}-1\right)}.
\end{equation}
By the above choice of $r$ and given $\Vert f \Vert_{D(A)} \leq r$ with $w \in \overline{B_{r}(0)}$, we have from \eqref{linear_Schrodinger_unitary_group_estimate} and \eqref{check_mapping_condition} that
\begin{align}
\begin{split}
\left\Vert S\left(-q(x)N(v + w,\overline{v + w})\right) \right\Vert_{C^{0}_{D(A)}}
\\ \leq
C_{0} C_{1} T K^{*}\Vert q\Vert_{H^{2}(\Omega)} \left(C_{1}\Vert f \Vert_{D(A)} + \Vert w \Vert_{C^{0}_{D(A)}}\right)^{m_0 + n_0} \\ \leq C_{0} C_{1} TK^{*} \Vert q \Vert_{H^{2}(\Omega)} \cdot \left(1 + C_{1}\right)^{m_0 + n_0} r^{m_{0} + n_{0}}\leq r,
\end{split}
\end{align}
which verifies \eqref{mapping_condition_0}. Observe that \eqref{stability_constants_nonlinear_well-posedness_0} relies on the assumption that $m_{0} + n_{0} > 1$.  \\ \\
\textit{Step 3 (Obtaining a Contraction Map).}\\
Now, we check that if $\Vert f \Vert_{D(A)} \leq r$ and $w \in \overline{B_{r}(0)}$ for $r > 0$ sufficiently small, then the map defined by 
\begin{equation} \label{contraction_candidate_definition_0}
w \longmapsto F(w,f) = S\left(-q(x)N(v + w,\overline{v + w})\right)
\end{equation}
is a contraction map on $\overline{B_{r}(0)}$. We find that
\begin{equation} \label{contraction_verification}
\begin{split}
\Vert F(w_1,f) - F(w_2,f)\Vert_{C^{0}_{D(A)}} \\
\leq C_{1} T K^{*}\Vert q \Vert_{H^{2}(\Omega)} \left\Vert N(v + w_1,\overline{v + w_1}) - N(v + w_2,\overline{v + w_2})\right\Vert_{C^{0}_{D(A)}},
\end{split}
\end{equation}
 where the above inequality was obtained by estimating as in \eqref{check_mapping_condition} through the unitary group estimate in \eqref{linear_Schrodinger_unitary_group_estimate} and the Banach algebra property of $D(A)$.\\ \indent
 In order to estimate the norm of
 \begin{equation} \label{contraction_verification_1}
 N(v + w_1,\overline{v + w_1}) - N(v + w_2,\overline{v + w_2}),
 \end{equation}
 we carry out an elementary calculation based on substituting the power series from \eqref{nonlinear_term_power_series}, which holds for $N(v + w_1,\overline{v + w_1})$ and $N(v + w_2,\overline{v + w_2})$ when $r > 0$ is sufficiently small. Adding and subtracting ``$(v+w_{1})^{l-m}(\overline{v+w_{2}})^{m}$" while fixing some $t \in [0, T]$, we see that the partial sum of \eqref{contraction_verification_1} for $k \leq l \leq K$ is bounded above in norm by
\begin{equation} \label{contraction_verification_2}
\begin{split}
\sum_{l = k}^{K}\Bigg(\sum_{m = 0}^{l}\frac{1}{m!(l-m)!}\left|\partial^{l-m}_{z_{1}}\partial^{m}_{z_{2}}N(0,0)\right|\Bigg[\bigg\Vert (v+w_{1})^{l-m}(\overline{v+w_{1}})^{m} - (v+w_{1})^{l-m}(\overline{v+w_{2}})^{m}  \bigg\Vert_{D(A)}\\
+\;\bigg\Vert (v+w_{1})^{l-m}(\overline{v+w_{2}})^{m} - (v+w_{2})^{l-m}(\overline{v+w_{2}})^{m}  \bigg\Vert_{D(A)} \Bigg]\Bigg).
\end{split}
\end{equation}
Moreover, we have from the Banach algebra property in \eqref{Banach_algebra_constant_0} that
\begin{equation} \label{contraction_verification_2_0_0}
\begin{split}
\sum_{l = k}^{K}\Bigg(\sum_{m = 0}^{l}\frac{1}{m!(l-m)!}\left|\partial^{l-m}_{z_{1}}\partial^{m}_{z_{2}}N(0,0)\right|\bigg\Vert (v+w_{1})^{l-m}(\overline{v+w_{1}})^{m} - (v+w_{1})^{l-m}(\overline{v+w_{2}})^{m}  \bigg\Vert_{D(A)} \Bigg)\\
\leq \sum_{l = k}^{K}\Bigg(\sum_{m = 0}^{l}\frac{1}{m!(l-m)!}\left|\partial^{l-m}_{z_{1}}\partial^{m}_{z_{2}}N(0,0)\right| \cdot (K^{*})^{l-m}\left\Vert v+w_{1} \right\Vert^{l-m}_{D(A)} \left\Vert (\overline{v + w_{1}})^{m} - (\overline{v + w_{2}})^{m} \right\Vert_{D(A)} \Bigg) \\
\leq \sum_{l = k}^{K}\Bigg(\sum_{m = 1}^{l}\frac{1}{m!(l-m)!}\left|\partial^{l-m}_{z_{1}}\partial^{m}_{z_{2}}N(0,0)\right| \cdot (K^{*})^{l-1}\left\Vert v+w_{1} \right\Vert^{l-m}_{D(A)} \Bigg(\sum_{j_{1} = 0}^{m - 1}\Vert v + w_{1} \Vert^{m - j_{1} - 1}_{D(A)} \Vert v + w_{2} \Vert^{j_{1}}_{D(A)}\Bigg) \Bigg)\\ \cdot \Vert w_{1} - w_{2} \Vert_{D(A)}.
\end{split}
\end{equation}
Since $\Vert v \Vert_{C^{0}_{D(A)}} \leq C_{1}r$,$\Vert w_{1} \Vert_{C^{0}_{D(A)}} \leq r$, $\Vert w_{2} \Vert_{C^{0}_{D(A)}} \leq r$, and $ 1 \leq m \leq l$, we see that \eqref{contraction_verification_2_0_0} is bounded above by
\begin{equation} \label{contraction_verification_2_0}
\sum_{l = k}^{K}\Bigg(\sum_{m = 0}^{l}\frac{1}{m!(l-m)!}\left|\partial^{l-m}_{z_{1}}\partial^{m}_{z_{2}}N(0,0)\right|\cdot l\left[(1+C_{1})K^{*}r \right]^{l-1} \Bigg) \cdot \Vert w_{1} - w_{2} \Vert_{D(A)}.
\end{equation}
We similarly have that the second term in \eqref{contraction_verification_2} is bounded above by \eqref{contraction_verification_2_0}.\\ \indent
Next, recall from Section \ref{Introduction} that we assumed there exists $\delta > 0$ such that $N = N(z_{1},z_{2})$ is analytic on the disk $B_{\delta}(0,0)$ defined by \eqref{nonlinear_term_disk_of_analyticity}. It follows from this that when $z_{1} = z_{2} = z$ is on the real axis, the series
\begin{equation} \label{convergence_of_derivative_series}
\sum_{l = k}^{\infty}\left(\sum_{m = 0}^{l} \frac{1}{m!(l-m)!} \left(\partial_{z_1}^{l-m}\partial_{z_2}^{m}N(0,0) \right) \cdot lz^{l-1} \right)
\end{equation}
for the real derivative converges absolutely. Combining \eqref{contraction_verification_2_0} - \eqref{convergence_of_derivative_series} while taking $K \rightarrow \infty$ and the supremum over $[0,T]$, it follows that for 
\begin{equation}
0 < r < \frac{\sqrt{\delta}}{\sqrt{2}(1 + C_{1})K^{*}}
\end{equation}
sufficiently small, the map defined by \eqref{contraction_candidate_definition_0} is a contraction map on $\overline{B_{r}(0)}$. We conclude as a result that the Contraction Mapping Principle yields a unique solution of the form $u = v + w$ to \eqref{quadraticStartEqn} in $\overline{B_{r}(0)}$.\\ \\
\textit{Step 4 (Stability Estimate).}\\
In the following, we estimate the solution by the initial condition using \eqref{check_mapping_condition} and by picking
\begin{equation}
0 < r \leq r_{m_{0},n_{0}} \cdot \left(\frac{1}{2C_{1}} \right)^{1/\left(m_{0} + n_{0} - 1\right)},
\end{equation}
where $r_{m_{0}, n_{0}}$ is defined in \eqref{stability_constants_nonlinear_well-posedness_0} and $C_{1} > 1$ comes from \eqref{linear_Schrodinger_unitary_group_estimate}.
\begin{equation} \label{initial_condition_bound}
\begin{split}
\Vert w \Vert_{C^{0}_{D(A)}} \leq C_{0}  C_{1}  T K^{*}\Vert q\Vert_{H^{2}(\Omega)} \left(\Vert v \Vert_{C^{0}_{D(A)}} + \Vert w \Vert_{C^{0}_{D(A)}}\right)^{m_0 + n_0} \\ \leq 
C_{0}  C_{1}  T K^{*}\Vert q\Vert_{H^{2}(\Omega)} \left( C_{1}\Vert f \Vert_{D(A)} + \Vert w \Vert_{C^{0}_{D(A)}} \right)   \left( C_{1}\Vert f \Vert_{D(A)}  + \Vert w \Vert_{C^{0}_{D(A)}}\right)^{m_0 + n_0 - 1}\\
\leq C_{0}  C_{1}  T K^{*}\Vert q\Vert_{H^{2}(\Omega)} \left(1 + C_{1}\right)^{m_0 + n_0 - 1}  r^{m_{0} + n_{0} - 1}  \left( C_{1}\Vert f \Vert_{D(A)}  + \Vert w \Vert_{C^{0}_{D(A)}}\right) \\
\leq C_{0}  C_{1}  T K^{*}\Vert q\Vert_{H^{2}(\Omega)}  \left(1 + C_{1}\right)^{m_0 + n_0}  r^{m_{0} + n_{0} - 1}  \left( C_{1}\Vert f \Vert_{D(A)}  + C_{1}\Vert w \Vert_{C^{0}_{D(A)}}\right)
\\
\leq \frac{1}{2C_{1}}\left( C_{1}\Vert f \Vert_{D(A)}  + C_{1}\Vert w \Vert_{C^{0}_{D(A)}}\right).
\end{split}
\end{equation}
Now, \eqref{initial_condition_bound} yields $\Vert w \Vert_{C^{0}_{D(A)}} \leq \Vert f \Vert_{D(A)}$ and consequently 
\begin{equation} \label{nonlinear_initial_condition_stability}
\Vert u \Vert_{C^{0}_{D(A)}} = \Vert v + w \Vert_{C^{0}_{D(A)}} \leq \left( 1 + C_{1}\right)\Vert f \Vert_{D(A)}.
\end{equation}
Note that by an argument similar to that of Theorem \ref{linear_well-posedness}, we have $u(x,t) \in C^{1}([0,T];L^{2}(\Omega))$.

\section{Tools for the Inverse Problem} \label{main_theorem_proof_tools}
\subsection{Linearization of the Schr\"odinger Equation} \label{rigorous_linearization}
In the rest of the paper, ``$\partial u/\partial \epsilon$" is taken to mean the derivative of the map from $\mathbb{R} \longrightarrow C^{0}_{D(A)}$ defined by
\begin{equation}\label{epsilon_derivative_meaning}
\epsilon \longmapsto u = u(x, t; \epsilon f),
\end{equation}
which denotes the unique solution to \eqref{quadraticStartEqn} with initial data $\epsilon f(x)$. The restriction to $\epsilon \in \mathbb{R}$ ensures that
\begin{equation*}
\frac{\partial \overline{u}}{\partial \epsilon} = \overline{\frac{\partial u}{\partial \epsilon}}.
\end{equation*}
For $l \in \mathbb{N}$ and $j = 1,2$, define in relation to \eqref{quadraticStartEqn_0}
\begin{equation} \label{definition_linearized_solution}
^{(l)}u_{\epsilon,j}(x,t) = \frac{\partial^{l} u_{\epsilon,j}}{\partial \epsilon^{l}}(x,t)\bigg\vert_{\epsilon = 0}.
\end{equation}
The left superscript indicates the order of the linearization with respect to $\epsilon \in \mathbb{R}$. \\ \indent
The following result establishes the differentiability of $u_{\epsilon,j}(x,t)$ satisfying \eqref{quadraticStartEqn_0} with respect to $\epsilon$. Its proof relies on showing the convergence of certain difference quotients. We refer the reader to the proof of Proposition A.4 in \cite{GPEquation_Linearization} where a similar idea is used.
\begin{lemma} \label{GPEqn_linearization_justification}
Let $N(z_1,z_2)$ be an analytic mapping near $(0,0)$ that satisfies the conditions in \eqref{nonlinearity_condition_0}-\eqref{nonlinearity_conditions} for some integer $k > 1$. If $u_{\epsilon,j}(x,t) \in C([0,T];D(A)) \cap C^1([0,T];L^2(\Omega))$ satisfies \eqref{quadraticStartEqn_0}, ${^{(1)}}u_{\epsilon,j}(x,t)$ exists in the sense of \eqref{epsilon_derivative_meaning} and satisfies
\begin{equation} \label{FirstLinear_for_p1}
\begin{cases} \left(i\partial_t + \Delta  + p_j(x) \right){^{(1)}}u_{\epsilon,j}(x,t) = 0 \text{ on } \Omega \times (0,T),\\ {^{(1)}}u_{\epsilon,j}(x,t) = 0 \text{ on } \partial \Omega \times (0,T),\\ {^{(1)}}u_{\epsilon,j}(x,0) = f(x) \text{ on } \Omega \times \{0\}.\end{cases}
\end{equation}
Additionally, ${^{(k)}}u_{\epsilon,j}(x,t)$ exists and satisfies\footnote{Recall that for a nonnegative pair of integers $a \geq b$, we define $\binom{a}{b} = a!/\left[(a-b)!\,b!\right]$.}
\begin{equation} \label{Second_Order}
\begin{cases}
\left(i\partial_t + \Delta + p_j(x)\right){^{(k)}}u_{\epsilon,j}(x,t) = -\sum_{m = 0}^{k}\binom{k}{m} \cdot q_j(x) \cdot \\
\left[\partial^{m}_{z_1}\partial^{k - m}_{z_2}N(0,0)\left({^{(1)}}u_{\epsilon,j}(x,t)\right)^{m}\left(\overline{{^{(1)}}u_{\epsilon,j}(x,t)}\right)^{k - m}\right] \text{ on } \Omega \times (0,T),\\ {^{(k)}}u_{\epsilon,j}(x,t) = 0 \text{ on } \partial \Omega \times (0,T),\\ {^{(k)}}u_{\epsilon,j}(x,0) = 0 \text{ on } \Omega \times \{0\}.
\end{cases}
\end{equation}
\end{lemma}
\begin{proof}
For $g \in C^{l}(\mathbb{R})$ with $l \in \mathbb{N}$, the following difference quotient formula can be established by induction:
\begin{equation} \label{order_l_difference_quotient}
\frac{\partial^{l}g}{\partial x^{l}} = \lim_{\eta \rightarrow 0} \frac{1}{{\eta}^{l}}\sum_{n = 0}^{l} \binom{l}{n} (-1)^{n}g(x+(l-n)\eta).
\end{equation}
Observe from Theorem \ref{linear_well-posedness} that the unique solution to \eqref{FirstLinear_for_p1} is $e^{tA}f(x)$. To establish \eqref{FirstLinear_for_p1}, we apply Duhamel's Formula in \eqref{var_param_formula} to \eqref{quadraticStartEqn_0} and observe the estimate
\begin{equation} \label{FirstLinear_for_p1_proof}
\begin{split}
\left \Vert \frac{u_{\epsilon,j}}{\epsilon} - e^{tA}f\right \Vert_{C^{0}_{D(A)}} = \left \Vert \frac{1}{\epsilon}\int_{0}^{t}e^{(t-s)A}\left[q_jN\left(u_{\epsilon,j},\overline{u_{\epsilon,j}}\right)(\cdot\,,s)\right] \,ds \right\Vert_{C^{0}_{D(A)}} \\ \leq  \frac{1}{\vert \epsilon \vert}  C_{0} C_{1} TK^{*}\Vert q_j\Vert_{H^{2}(\Omega)}\Vert u_{\epsilon,j} \Vert^{m_0 + n_0}_{C^{0}_{D(A)}}
\\ \leq 
\frac{\left|(1 + C_{1})\epsilon\right|^{m_0 + n_0}}{\vert \epsilon \vert}  C_{0}  C_{1}  TK^{*}\Vert q_j\Vert_{H^{2}(\Omega)}\Vert f\Vert^{m_0 + n_0}_{D(A)}.
\end{split}
\end{equation}
The above inequalities use \eqref{nonlinearity_conditions}, the estimate in \eqref{linear_Schrodinger_unitary_group_estimate}, and the fact that \eqref{nonlinear_initial_condition_stability} holds for $\epsilon$ sufficiently small. Taking $\epsilon \rightarrow 0$, we obtain \eqref{FirstLinear_for_p1} upon noting that $m_{0} + n_{0} > 1$. It remains to establish \eqref{Second_Order}. \\ \indent
The well-posedness of \eqref{Second_Order} follows from Theorem \ref{linear_well-posedness} and \eqref{var_param_formula}. Applying \eqref{var_param_formula} to \eqref{Second_Order} and recalling the difference quotient formula \eqref{order_l_difference_quotient}, we consider
\begin{multline} \label{convergence_requirement_GPEqn_linearization}
\Bigg \Vert \frac{1}{\epsilon^{k}}\sum_{n = 0}^{k} \binom{k}{n} (-1)^{n}u_{(k - n)\epsilon,j} \\
-\; i\int_{0}^{t} e^{(t-s)A}\Bigg[q_j\sum_{m = 0}^{k}\binom{k}{m} \partial^{m}_{z_1}\partial^{k - m}_{z_2}N(0,0)\left({^{(1)}}u_{\epsilon,j}\right)^{m}\left(\overline{{^{(1)}}u_{\epsilon,j}}\right)^{k - m}(\cdot\,,s) \Bigg]\, ds
\Bigg \Vert_{C^{0}_{D(A)}}
\end{multline}
as $\epsilon \rightarrow 0$.
By \eqref{var_param_formula}, we have that
\begin{equation} \label{var_param_formula_difference_quotient}
u_{(k - n)\epsilon,j}(x,t) = (k - n)\epsilon e^{tA}f(x) \;+\; i\int_{0}^{t} e^{(t-s)A}\left[q_jN\left(u_{(k - n)\epsilon,j}, \overline{u_{(k - n)\epsilon,j}}\right)(x\,,s)\right]\,ds.
\end{equation}
When substituting \eqref{var_param_formula_difference_quotient} into \eqref{convergence_requirement_GPEqn_linearization}, observe that
\begin{equation} \label{combinatorial_identity}
\sum_{n = 0}^{k} \binom{k}{n}(-1)^{n}(k - n)e^{tA}f(x) 
= \left[\sum_{n = 0}^{k - 1} k\binom{k - 1}{n}(-1)^{n}\right]e^{tA}f(x) = 0
\end{equation}
as a result of the identities
\begin{equation}
(k - n)\binom{k}{n} = k\binom{k - 1}{n}
\quad \text{ and } \quad
\sum_{n = 0}^{k - 1} \binom{k - 1}{n}(-1)^{n} = 0.
\end{equation}
It remains to consider
\begin{equation} \label{GPEqn_Step_1}
\begin{split}
\Bigg \Vert \frac{1}{\epsilon^{k}}\int_{0}^{t} e^{(t-s)A}\Bigg[q_j\sum_{n=0}^{k} \binom{k}{n}(-1)^{n}N\left(u_{(k - n)\epsilon,j}, \overline{u_{(k - n)\epsilon,j}}\right)(\cdot\,,s) \Bigg]\,ds\\
- \int_{0}^{t} e^{(t-s)A}\Bigg[q_j(x)\sum_{m = 0}^{k}\binom{k}{m} \partial^{m}_{z_1}\partial^{k - m}_{z_2}N(0,0)\left({^{(1)}}u_{\epsilon,j}\right)^{m}\left(\overline{{^{(1)}}u_{\epsilon,j}}\right)^{k - m}(\cdot\,,s) \Bigg]\, ds \Bigg \Vert_{C^{0}_{D(A)}}
\end{split}
\end{equation}
as $\epsilon \rightarrow 0$. For $\epsilon$ sufficiently small, we can substitute the power series in \eqref{nonlinear_term_power_series} for $N\left(u_{(k - n)\epsilon,j},\overline{u_{(k - n)\epsilon,j}}\right)$ in \eqref{GPEqn_Step_1}. We can then bound \eqref{GPEqn_Step_1} above by $I_{1} + I_{2}$ where
\begin{equation} \label{GPEqn_Step_2}
\begin{split}
I_1 = \Bigg \Vert 
\frac{1}{\epsilon^{k}}\int_{0}^{t} e^{(t-s)A}\Bigg[q_{j}\sum_{n=0}^{k} \binom{k}{n}(-1)^{n}\Bigg(\sum_{h = 0}^{k} \frac{1}{h!(k - h)!}\left(\partial^{k -h}_{z_1}\partial^{h}_{z_2}N(0,0)\right) \cdot
u^{k - h}_{(k - n)\epsilon,j}\overline{u_{(k - n)\epsilon,j}}^{h} \Bigg)(\cdot\,,s)\Bigg]\,ds \\ -\;
\int_{0}^{t}e^{(t-s)A}\Bigg[q_{j}\sum_{m = 0}^{k}\binom{k}{m} \partial^{m}_{z_1}\partial^{k - m}_{z_2}N(0,0)\left({^{(1)}}u_{\epsilon,j}\right)^{m}\left(\overline{{^{(1)}}u_{\epsilon,j}}\right)^{k - m}(\cdot\,,s)\Bigg] \,ds
\Bigg \Vert_{C^{0}_{D(A)}}
\end{split}
\end{equation}
and
\begin{equation}
\begin{split}
I_{2} = \Bigg \Vert \frac{1}{\epsilon^{k}}\int_{0}^{t} e^{(t-s)A}
\Bigg[q_{j}\sum_{n = 0}^{k}\binom{k}{n} (-1)^{n} \Bigg(\sum_{l = k + 1}^{\infty} \sum_{h = 0}^{l} \frac{1}{h!(l-h)!}\left(\partial^{l-h}_{z_1}\partial^{h}_{z_2}N(0,0)\right) \cdot \\ u^{l-h}_{(k - n)\epsilon,j}(x,s)\overline{u_{(k - n)\epsilon,j}}^{h}\Bigg)(\cdot\,,s) \Bigg]
\,ds\Bigg \Vert_{C^{0}_{D(A)}},
\end{split}
\end{equation}
where we have grouped the $l = k$ term of the power series with $I_1$ and the terms for $l \geq k + 1$ with $I_2$.\\ \indent
Let us write
\begin{equation} \label{splitting_absolute_convergence}
\epsilon f(x) = \epsilon_{1} \cdot \epsilon_{2}f(x)
\end{equation}
where $\epsilon_{1} \rightarrow 0$ and with a fixed $\epsilon_{2} \in \mathbb{R}$ such that $N(\epsilon_{2}f, \epsilon_{2}f)$ converges absolutely with respect to $\Vert \cdot \Vert_{C^{0}_{D(A)} \times C^{0}_{D(A)}}$ upon regarding $f(x) = f(x,\cdot) \in C^{0}_{D(A)}$. Consequently, observe that when $0 < |\epsilon_{1}| < 1/[K^{*}\left(1 + C_{1}\right)k]$, we have
\begin{multline} \label{dominated_convergence_0}
%\begin{split}
\frac{1}{|\epsilon|^{k}}\Bigg(\sum_{l = k + 1}^{\infty} \sum_{h=0}^{l}   \frac{1}{h!(l-h)!}\left|\partial^{l-h}_{z_1}\partial^{h}_{z_2}N(0,0)\right|  %\\
(K^{*})^{l-1}\left(1 + C_{1} \right)^{l}(k - n)^{l} \Vert \epsilon f \Vert^{l}_{D(A)} \Bigg) \\=
\frac{1}{|\epsilon_{2}|^{k}}\Bigg(\sum_{l = k + 1}^{\infty} \sum_{h=0}^{l} \frac{1}{h!(l-h)!}\left|\partial^{l-h}_{z_1}\partial^{h}_{z_2}N(0,0)\right|  %\\
(K^{*})^{l-1}\left(1 + C_{1} \right)^{l}(k - n)^{l}|\epsilon_{1}|^{l-k} \Vert \epsilon_{2}f \Vert^{l}_{D(A)} \Bigg)
< \infty.
%\end{split}
\end{multline}
To see this, note that for $0 < |\epsilon_{1}| < 1/[K^{*}\left(1 + C_{1}\right)k]$, $0 \leq n \leq k$, and $l \geq k + 1$,
\begin{equation} \label{dominated_convergence_technical_point}
(K^{*})^{l-1}\left(1 + C_{1} \right)^{l}(k - n)^{l}|\epsilon_{1}|^{l - k} < \left[\left(1 + C_{1} \right)k\right]^{k}(K^{*})^{k-1}. 
\end{equation}
Using \eqref{splitting_absolute_convergence} along with the estimates in \eqref{linear_Schrodinger_unitary_group_estimate} and \eqref{nonlinear_initial_condition_stability} (which holds for $\epsilon$ sufficiently small), we have the estimate
\begin{align} \label{GPEqn_Step_3}
\begin{split}
I_2 \leq \frac{1}{|\epsilon_{2}|^{k}} C_{1}  T K^{*}\Vert q_j\Vert_{H^{2}(\Omega)}  \sum_{n = 0}^{k}\binom{k}{n}\Bigg(\sum_{l = k + 1}^{\infty} \sum_{h=0}^{l}   \frac{1}{h!(l-h)!}\left|\partial^{l-h}_{z_1}\partial^{h}_{z_2}N(0,0)\right| \cdot \\
(K^{*})^{l-1}\left(1 + C_{1} \right)^{l}(k - n)^{l}|\epsilon_{1}|^{l-k} \Vert \epsilon_{2} f \Vert^{l}_{D(A)} \Bigg)
\end{split}
\end{align}
Taking $0 < |\epsilon_{1}| < 1/[K^{*}\left(1 + C_{1} \right)k]$ and recalling \eqref{dominated_convergence_0}, the Lebesgue Dominated Convergence Theorem shows that $I_2 \rightarrow 0$ as $\epsilon_{1} \rightarrow 0$. \\ \indent
Let's now investigate $I_1$ as $\epsilon \rightarrow 0$. Using an argument similar to that of \eqref{FirstLinear_for_p1_proof}, observe for $\alpha \in \mathbb{C}$ that
\begin{equation} \label{GP_Eqn_Step_3.5}
\frac{u_{\alpha\epsilon,j}}{\epsilon} \longrightarrow \alpha \cdot {^{(1)}}u_{\epsilon,j}
\end{equation}
as $\epsilon \rightarrow 0$. We see then that
\begin{align} \label{GPEqn_Step_4}
\begin{split}
\frac{1}{\epsilon^{k}}u^{k -h}_{(k - n)\epsilon,j}\overline{u_{(k - n)\epsilon,j}}^{h} = \left(\frac{u_{(k - n)\epsilon,j}}{\epsilon}\right)^{k - h}\left(\overline{\frac{u_{(k - n)\epsilon,j}}{\epsilon}}\right)^{h} \\ \longrightarrow (k - n)^{k} \cdot 
{^{(1)}}u_{\epsilon,j}^{k - h} \cdot \overline{^{(1)}u_{\epsilon,j}}^{h}
\end{split}
\end{align}
as $\epsilon \rightarrow 0$. Note also the identity (e.g. see Theorem 3.3.1 and Corollary 3.3.2 in \cite{Combinatorics_Reference})
\begin{equation} \label{GPEqn_Step_5}
\sum_{n=0}^{k}\binom{k}{n}(-1)^{n}(k - n)^{k} = k!
\end{equation}
It follows from \eqref{GPEqn_Step_4} and \eqref{GPEqn_Step_5} that $I_1 \rightarrow 0$ as $\epsilon \rightarrow 0$, which establishes \eqref{Second_Order}.
\end{proof}
\begin{remark}
Using the argument of Lemma \ref{GPEqn_linearization_justification}, one can show for $2 \leq l \leq k - 1$ that ${^{(l)}}u_{\epsilon,j} \equiv 0$.
\end{remark}

\subsection{Carleman Estimates for the Linear Schr\"odinger Operator} \label{main_previous_literature_Carleman_estimate}
\begin{figure}[h]
    \centering
    \includegraphics[width=0.35\linewidth]{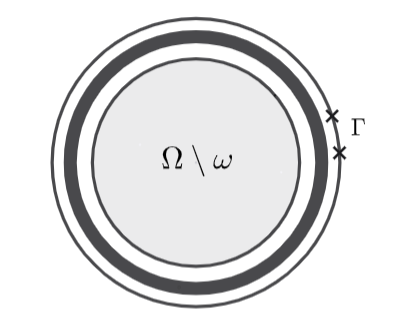}
    \caption{The black region within $\omega$ schematically represents $\overline{\Omega_{2} \setminus \Omega_{1}}$, which appears in \eqref{logarithmic_observability_inequality_statement}.}
    \label{Figure_1}
\end{figure}

In the rest of the paper, ``RHS" abbreviates ``right-hand side" and ``LHS" abbreviates ``left-hand side". To prove Theorems \ref{main_theorem}-\ref{main_theorem_partial_data}, we shall rely on the Carleman estimates established in \cite{StartingCarlemanEstimate,MainMethod}.\\ \indent
We begin by noting the following estimate for $L = i\partial_{t} + \Delta + p(x)$ from Proposition 1 in \cite{StartingCarlemanEstimate}:
\begin{proposition} \label{standard_Schrodinger_Carleman_estimate_1}
Pick some $x_0 \notin \overline{\Omega}$, $T_{1} \in (0, T]$, and consider $\psi(x) = |x-x_0|^2$. Define the functions 
\begin{align} \label{weight_function}
\theta(x,t) = \frac{e^{\lambda \psi(x)}}{T_{1}^{2}-t^{2}} \quad \text{and} \quad \phi(x,t) = \frac{e^{2\lambda\Vert \psi\Vert_{\infty}} - e^{\lambda \psi(x)}}{T_{1}^{2}-t^{2}}
\end{align}
on $\Omega \times (-T_{1},T_{1})$, where $\Vert \psi \Vert_{\infty} = \Vert \psi \Vert_{L^{\infty}(\Omega)}$. Suppose $u \in L^2(-T_{1},T_{1};H^1_{0}(\Omega))$, $Lu \in L^2(\Omega \times (-T_{1},T_{1}))$, and $\partial u/\partial \nu \in L^2(-T_{1},T_{1};L^2(\partial \Omega))$, and that there exists $M > 0$ such that $\Vert p\Vert_{L^{\infty}(\Omega)} \leq M$. We denote $ \Gamma_{1} = \{x \in \partial\Omega \,\big|\, (\nabla \psi \cdot \nu)(x) \geq 0\}$, $R_{1} = e^{s\phi}\left(i\partial_t + \Delta + s^2|\nabla \phi|^2\right)e^{-s\phi}$, and $R_{2} = e^{s\phi}\left(is(\partial_t \phi)+2s(\nabla \phi \,\cdot \nabla) + s\Delta \phi\right)e^{-s\phi}$. \\ \indent Consider $\Gamma_{0} \subseteq \partial \Omega$ such that $\Gamma_{1} \subseteq \Gamma_{0}$. Then, there exists $s_0 = s_{0}(\Omega, T_{1}) > 0$, $\Lambda_0 = \Lambda_{0}(\Omega, T_{1}) > 0$, and some constant $C = C(\Omega, \Gamma_{0}, T_{1}, M, s_{0}, \Lambda_{0}) > 0$ independent of $s\in [s_0, \infty)$ and $\lambda \in [\Lambda_0, \infty)$ so that
\begin{multline} \label{OldCarlemanEstimate_1}
s^3 \lambda^4\Vert e^{-s\phi}u\Vert_{L^2(\Omega \times (-T_{1},T_{1}))}^2 + s\lambda\Vert e^{-s\phi}\nabla u\Vert_{L^2(\Omega \times (-T_{1},T_{1}))}^2 + \Vert e^{-s\phi}R_1(u) \Vert_{L^2(\Omega \times (-T_{1},T_{1}))}^2 \\ +\; \Vert e^{-s\phi}R_2(u)\Vert_{L^2(\Omega \times (-T_{1},T_{1}))}^2 \leq C\bigg(\Vert e^{-s\phi} Lu\Vert_{L^2(\Omega \times (-T_{1},T_{1}))}^2\;+\;  s\lambda\left\Vert \sqrt{\theta}e^{-s\phi}\frac{\partial u}{\partial \nu} \sqrt{\nabla \psi \cdot \nu}\right\Vert_{L^2(\Gamma_{0} \times (-T_{1},T_{1}))}^2\bigg).
\end{multline}
\end{proposition}
Observe that $\vert \nabla \psi \cdot \nu \vert$ is uniformly bounded on $\overline{\Omega}$. Also, note that for $\Lambda_0(\Omega, T_{1}) > 0$ sufficiently large and $\lambda \geq \Lambda_{0}$, we have
\begin{equation} \label{weight_function_bound_0}
\frac{e^{2\lambda \Vert \psi \Vert_{\infty}}}{T_{1}^2 - t^2} \geq \frac{\left(e^{\lambda \psi(x)}\right)^2}{T_{1}^2 - t^2} \geq \frac{2e^{\lambda \psi(x)}}{T_{1}^2 - t^2}.
\end{equation}
Then, we have for $\lambda \in [\Lambda_0, \infty)$ that
\begin{equation} \label{weight_function_bound_1}
\phi(x,t) \geq \frac{2e^{\lambda \psi(x)} - e^{\lambda \psi(x)}}{T_{1}^2 - t^2} = \theta(x,t).
\end{equation}
As a result, we can majorize the boundary term in \eqref{OldCarlemanEstimate_1} by
\begin{equation} \label{boundary_term_bound}
\left\Vert \sqrt{\theta}e^{-s\phi}\frac{\partial u}{\partial \nu} \sqrt{\nabla \psi \cdot \nu}\right\Vert_{L^2(\Gamma_{0} \times (-T_{1},T_{1}))}^2 \leq C \left\Vert \frac{\partial u}{\partial \nu} \right\Vert_{L^2(\Gamma_{0} \times (-T_{1},T_{1}))}^2
\end{equation}
for some constant $C = C(\Omega, \Gamma_{0}, T_{1}, M, s_{0},\Lambda_0) > 0$.\\
\indent We now state another Carleman estimate for $L = i\partial_{t} + \Delta + p(x)$ proved in \cite{MainMethod}. We begin by defining neighborhoods $\{\omega_{1}, \omega_{2}\}$ of $\partial \Omega$ such that
\begin{equation}
\overline{\omega_{j + 1}} \subset \omega_{j}
\end{equation}
for $j = 0,1$ (where $\omega_{0} = \omega$ is introduced in the assumptions of Theorem \ref{main_theorem_partial_data}). We also define
\begin{equation} \label{partial_data_proof_subset_1}
\Omega_{j} = \Omega \setminus \overline{\omega_{j}} \quad (j = 1,2).
\end{equation}
We refer the reader to Figure \ref{Figure_1} for an illustration.\\ \indent
The following is from Proposition $2.1$ in \cite{MainMethod}.
\begin{proposition} \label{standard_Schrodinger_Carleman_estimate_0}
Let $p(x) \in L^{\infty}(\Omega)$ with $\Vert p\Vert_{L^{\infty}(\Omega)} \leq M$ for some $M > 0$. Picking a point $x_{0} \notin \overline{\Omega}$ and $T_{1} \in (0, T]$, define $\psi(x)$ and $\phi(x,t)$ as in Proposition \ref{standard_Schrodinger_Carleman_estimate_1}. Now, there exist $\lambda_{*} = \lambda_{*}(\Omega, T_{1}, \Omega_{1}, \Omega_{2}) > 0,\, s_{*} = s_{*}(\Omega,T_{1},\Omega_{1}, \Omega_{2}) > 0$ such that for all $\lambda \geq \lambda_{*},\, s \geq s_{*},$ and a constant $C = C(\Omega, T_{1}, M,s_{*},\lambda_{*}) > 0$,
\begin{equation} \label{OldCarlemanEstimate_2}
\begin{split}
s^3 \lambda^4\Vert e^{-s\phi}u\Vert_{L^2(\Omega_{1} \times (-T_{1},T_{1}))}^2 + s\lambda\Vert e^{-s\phi}\nabla u\Vert_{L^2(\Omega_{1} \times (-T_{1},T_{1}))}^2
\leq
C\Vert e^{-s\phi} Lu\Vert_{L^2(\Omega \times (-T_{1},T_{1}))}^2 \\ 
+\; C\left(s^3 \lambda^4\Vert e^{-s\phi}u\Vert_{L^2(\Omega_{2} \setminus \Omega_{1} \times (-T_{1},T_{1}))}^2 + s\lambda\Vert e^{-s\phi}\nabla u\Vert_{L^2(\Omega_{2} \setminus \Omega_{1} \times (-T_{1},T_{1}))}^2 \right).
\end{split}
\end{equation}
The above holds for $u \in L^{2}\left(-T_{1},T_{1};\,H^{1}(\Omega) \right)$ such that $Lu \in L^{2}\left(\Omega \times (-T_{1},T_{1})\right)$.
\end{proposition}
The proof of \eqref{OldCarlemanEstimate_2} is based on the estimate
\begin{equation} \label{OldCarlemanEstimate_2_step_0}
s^3 \lambda^4\Vert e^{-s\phi}w\Vert_{L^2(\Omega \times (-T_{1},T_{1}))}^2 + s\lambda\Vert e^{-s\phi}\nabla w\Vert_{L^2(\Omega \times (-T_{1},T_{1}))}^2
 \leq C\Vert e^{-s\phi} Lw\Vert_{L^2(\Omega \times (-T_{1},T_{1}))}^2,
\end{equation}
which is a special case of \eqref{OldCarlemanEstimate_1} for $w(x,t) \in C^{\infty}(-T_{1}, T_{1}; C^{\infty}_{0}(\Omega))$. It also holds for $w(x,t) \in L^{2}\left(-T_{1},T_{1};\,H^{1}_{0}(\Omega) \right)$ by density in $L^{2}\left(-T_{1},T_{1};\,H^{1}(\Omega) \right)$. The estimate in \eqref{OldCarlemanEstimate_2} then follows as in \cite{MainMethod} by applying \eqref{OldCarlemanEstimate_2_step_0} to $w(x,t) = \chi(x)u(x,t)$ where $\chi(x)$ is a smooth function satisfying
\begin{equation}
\chi(x) = \begin{cases}1,\;x \in \overline{\Omega_{1}}\\
0,\;x \in \Omega \setminus \overline{\Omega_{2}}.
\end{cases}
\end{equation}
\begin{remark} \label{Carleman_weight_observation_at_zero_technical_point}
\hfill
\begin{itemize}
\item[(a)] Observe that since we have $\Omega_{2} = \Omega_{1} \cup \left(\Omega_{2}\setminus \Omega_{1} \right)$, the estimate in \eqref{OldCarlemanEstimate_2} also holds with ``$\Omega_{2} \times (-T_{1}, T_{1})$" on the LHS in place of ``$\Omega_{1} \times (-T_{1}, T_{1})$" for some constant $C = C(M,s_{*},\lambda_{*}) > 1$.
\item[(b)] Regarding Proposition \ref{main_previous_literature_Carleman_estimate}, the parameter dependence $s_{0} = s_{0}(\Omega, T_{1})$ and $\Lambda_{0} = \Lambda_{0}(\Omega, T_{1})$ is not explicitly stated in Proposition $1$ from \cite{StartingCarlemanEstimate}. However, the dependence is apparent from the proof provided in \cite{StartingCarlemanEstimate}, where the ability to absorb terms into $s^3 \lambda^4\Vert e^{-s\phi}u\Vert_{L^2(\Omega \times (-T_{1},T_{1}))}^2$ and $s \lambda\Vert e^{-s\phi}\nabla u\Vert_{L^2(\Omega \times (-T_{1},T_{1}))}^2$ depends on bounds on for the function $\theta$ (defined by \eqref{weight_function}) over $\Omega \times (-T_{1}, T_{1})$. 
\item[(c)] We comment on the reason for working on the symmetric interval $(-T_{1}, T_{1})$. In \eqref{weight_function}, it is possible to replace ``$T_{1}^{2}-t^{2}$" with ``$t(T_{1} - t)$" and follow the proof in \cite{StartingCarlemanEstimate} to obtain an estimate where the norms are taken on $\Omega \times (0, T_{1})$ and $\Gamma_{1} \times (0, T_{1})$. However, this choice of weight function does not satisfy $e^{-s\phi(x,0)} \neq 0$, which is necessary for the argument in Section \ref{main_theorem_proof}. In Section \ref{main_theorem_proof}, we will introduce a linearized Schr\"{o}dinger equation crucial to the proofs of Theorem \ref{main_theorem}-\ref{main_theorem_partial_data}. In order to work on the symmetric interval, we will extend the solution to this linearized equation in a way that enables us to conveniently control the Neumann data on $(-T_{1}, 0)$. We refer to Remark \ref{symmetric_extension_remark} for further details on the extension.
\end{itemize}
\end{remark}

\section{Proofs of Theorem \ref{main_theorem}-\ref{main_theorem_partial_data}}\label{main_theorem_proof}
We shall now prove the main results, initially focusing on \eqref{main_theorem_stable_determination} and \eqref{main_theorem_partial_stable_determination} dealing with the stable determination of $q_{j}(x)$. Before that we note some preliminary lemmas, beginning with an equation for which the initial data contains information about $q_{1}(x) - q_{2}(x)$.
\begin{remark} \label{notation_0}
In the rest of the paper, given a function $g(x,t)$, we adopt the convention that $n\cdot \left[ g(x,t)\right]^{n-1} = 0$ when $n = 0$.
\end{remark}
\begin{lemma} \label{Partial_t_Second_Order_Difference_lemma}
In relation to \eqref{quadraticStartEqn_0}, suppose that $p_{1}(x) = p_{2}(x) = p(x) \in H^{2}(\Omega)$ and $q_{j}(x) \in H^{2}(\Omega)\; (j = 1,2)$. Denoting 
\begin{equation} \label{Partial_t_Second_Order_Difference_source_term}
\begin{split}
v(x,t) = -\sum_{m = 0}^{k}\binom{k}{m} \cdot 
\partial^{m}_{z_1}\partial^{k - m}_{z_2}N(0,0) \cdot \\  \bigg[m \cdot \partial_{t} \left({^{(1)}}u_{\epsilon}(x,t)\right) \cdot \left({^{(1)}}u_{\epsilon}(x,t)\right)^{m-1}  \cdot  \left(\overline{^{(1)}u_{\epsilon}(x,t)}\right)^{k - m} + \\
(k - m) \cdot \partial_{t}\left(\overline{{^{(1)}}u_{\epsilon}(x,t)}\right) \cdot \left({^{(1)}}u_{\epsilon}(x,t)\right)^{m} \cdot 
\left(\overline{{^{(1)}}u_{\epsilon}(x,t)}\right)^{k - m - 1}
\bigg]
\end{split}
\end{equation}
and defining $r(x,t) = \partial_{t}\left({^{(k)}}u_{\epsilon,1} - {^{(k)}}u_{\epsilon,2}\right)(x,t)$, we have
\begin{align} \label{Partial_t_Second_Order_Difference}
\begin{cases}
\left(i\partial_t + \Delta + p(x)\right) r(x,t) = \left[q_1(x) - q_2(x)\right]v(x,t) \text{ on } \Omega \times (0,T),\\ r(x,t) = 0 \text{ on } \partial \Omega \times (0,T),\\ r(x,0) = i\sum_{m = 0}^{k}\binom{k}{m} \cdot \left[q_1(x) - q_2(x)\right] \cdot
\partial^{m}_{z_1}\partial^{k - m}_{z_2}N(0,0) \cdot \left[f(x)\right]^{k} \text{ on } \Omega \times \{0\}.
\end{cases}
\end{align}
Moreover, we have $r(x,t) \in C([0,T];D(A)) \cap C^{1}([0,T];L^2(\Omega))$ for $f(x) \in H^{4}(\Omega) \cap H^{3}_{0}(\Omega)$.
\end{lemma}
 To see that $r(x,t) \in C([0,T];D(A)) \cap C^{1}([0,T];L^2(\Omega))$ for $f(x) \in H^{4}(\Omega) \cap H^{3}_{0}(\Omega)$, first note that ${^{(1)}}\tilde{u_{\epsilon}}(x,t) = \partial_{t}\left({^{(1)}}u_{\epsilon}(x,t)\right)$ satisfies
\begin{gather} \label{partial_t_FirstLinear_for_p1}
\begin{cases}
\left(i\partial_t  + \Delta  + p(x)\right) {^{(1)}}\tilde{u_{\epsilon}}(x,t) = 0 \text{ on } \Omega \times (0,T),\\ {^{(1)}}\tilde{u_{\epsilon}}(x,t) = 0 \text{ on } \partial \Omega \times (0,T),\\ {^{(1)}}\tilde{u_{\epsilon}}(x,0)  = i\Delta f(x) + ip(x)f(x) \text{ on } \Omega \times \{0\}.
\end{cases}
\end{gather}
This follows by considering the weak formulation of \eqref{partial_t_FirstLinear_for_p1} from which we can show that ${^{(1)}}u_{\epsilon}(\cdot\,,t)$ and the mapping
\begin{equation}
t \longmapsto \int_{0}^{t} {^{(1)}}\tilde{u_{\epsilon}}(\cdot\,,s)\,ds
\end{equation}
satisfy the same equation. Now, taking $f(x) \in H^{4}(\Omega) \cap H^{3}_{0}(\Omega)$, we obtain ${^{(1)}}u_{\epsilon}(x,t) \in C^{1}([0,T];D(A)) \cap C^{2}([0,T];L^{2}(\Omega))$ upon applying Theorem \ref{linear_well-posedness} to \eqref{partial_t_FirstLinear_for_p1}. Now, applying \eqref{var_param_formula} to \eqref{Partial_t_Second_Order_Difference} yields $r(x,t) \in C([0,T];D(A)) \cap C^{1}([0,T];L^2(\Omega))$. The fact that $r(x,t) = \partial_{t}\left({^{(k)}}u_{\epsilon,1} - {^{(k)}}u_{\epsilon,2}\right)(x,t)$ can again be justified by considering the weak formulation of \eqref{Partial_t_Second_Order_Difference} and recalling \eqref{Second_Order}.\\ \indent
Now, we note some estimates for $r(x,t)$ and $v(x,t)$ that will be convenient in the later discussion.
\begin{lemma} \label{L2_H2_Partial_t_Second_Order_Difference_estimates}
For $j = 1,2$, consider $p_{j}(x),\,q_{j}(x) \in H^{2}(\Omega)$ with $\Vert p_{j} \Vert_{H^{2}(\Omega)} \leq M$ and $q_{1}(x) - q_{2}(x) \in H^{4}(\Omega)$. We also take $f(x) \in H^{4}(\Omega) \cap H^{3}_{0}(\Omega)$. For $r(x,t)$ a solution to \eqref{Partial_t_Second_Order_Difference} and $v(x,t)$ defined by \eqref{Partial_t_Second_Order_Difference_source_term}, we have for each fixed $t \in (0, T)$ the estimates
\begin{equation} \label{L2_H2_Partial_t_Second_Order_Difference_estimates_statement}
\begin{split}
\Vert v(\cdot\,,t) \Vert_{H^{2j}(\Omega)} \leq C(\Omega, M)\Vert f\Vert^{k}_{H^{4}(\Omega)}\: (j = 0,1),\\
\Vert r(\cdot\,,t) \Vert_{H^{2j}(\Omega)} \leq C(\Omega,M,T)\Vert q_{1}- q_{2} \Vert_{H^{2j}(\Omega)}\Vert f\Vert^{k}_{H^{4}(\Omega)}\: (j = 0,1),\\
\Vert \partial_{t}r(\cdot\,,t) \Vert_{H^{2j}(\Omega)} \leq C(\Omega, M)\Vert q_{1}- q_{2} \Vert_{H^{2j + 2}(\Omega)}\Vert f\Vert^{k}_{H^{4}(\Omega)}\: (j = 0, 1).
\end{split}
\end{equation}
\end{lemma}
\begin{proof}
Recall from Theorem \ref{linear_well-posedness} that when $g(x) \in D(A)$, we have the estimates
\begin{equation} \label{linear_well-posedness_result_estimates}
\begin{split}
\Vert e^{tA}g\Vert_{L^{2}(\Omega)} = \Vert g\Vert_{L^{2}(\Omega)},\\
\Vert Ae^{tA}g\Vert_{L^{2}(\Omega)} = \Vert Ag\Vert_{L^{2}(\Omega)},\\
\Vert e^{tA}g\Vert_{H^{2}(\Omega)} \leq C(\Omega,M)\Vert g\Vert_{H^{2}(\Omega)},
\end{split}
\end{equation}
which hold for each $t \in (0, T)$. Thus, we have from \eqref{FirstLinear_for_p1} and \eqref{partial_t_FirstLinear_for_p1} the estimates
\begin{equation}
\begin{split}
\left\Vert {^{(1)}}u_{\epsilon}(\cdot\,,t)  \right\Vert_{L^{2}(\Omega)} = \Vert f\Vert_{L^{2}(\Omega)},\\
\left\Vert {^{(1)}}u_{\epsilon}(\cdot\,,t)  \right\Vert_{H^{2}(\Omega)} \leq C(\Omega, M)\Vert f \Vert_{H^{2}(\Omega)},\\
\left\Vert \partial_{t}\left({^{(1)}}u_{\epsilon}(\cdot\,,t)\right) \right\Vert_{L^{2}(\Omega)} \leq C(M)\Vert f\Vert_{H^{2}(\Omega)},\\
\left\Vert  \partial_{t}\left({^{(1)}}u_{\epsilon}(\cdot\,,t) \right)\right\Vert_{H^{2}(\Omega)} \leq C(\Omega, M)\Vert f\Vert_{H^{4}(\Omega)}.
\end{split}
\end{equation}
Observe then from the above and the Sobolev Embedding Theorem that we have the estimates
\begin{equation}
\begin{split}
\left\Vert {^{(1)}}u_{\epsilon}(\cdot\,,t)  \right\Vert_{L^{\infty}(\Omega)} \leq C(\Omega, M)\Vert f\Vert_{H^{2}(\Omega)},\\
\left\Vert  \partial_{t}\left({^{(1)}}u_{\epsilon}(\cdot\,,t) \right)\right\Vert_{L^{\infty}(\Omega)} \leq C(\Omega, M)\Vert f \Vert_{H^{4}(\Omega)}.
\end{split}
\end{equation}
Now, from \eqref{Partial_t_Second_Order_Difference_source_term}, we obtain the estimate
\begin{equation} \label{Partial_t_Second_Order_Difference_source_term_L2_estimate}
\Vert v(\cdot\,,t) \Vert_{H^{2j}(\Omega)} \leq C(\Omega, M)\Vert f\Vert^{k}_{H^{4}(\Omega)}\: (j = 0,1),
\end{equation}
where the $\Vert v(\cdot\,,t) \Vert_{H^{2}(\Omega)}$ estimate uses the Banach algebra property of $H^{2}(\Omega)$.\\ \indent
Next, using the unitary group $\{e^{tA}\}_{t \in \mathbb{R}}$ from Theorem \ref{linear_well-posedness}, we can apply Duhamel's Formula to \eqref{Partial_t_Second_Order_Difference} to obtain
\begin{equation}
\begin{split}
r(x,t) = e^{tA}r(x,0) - i\int_{0}^{t}e^{(t-s)A}\left[q_1(x) - q_2(x)\right]v(x,s)\,ds,\\
\partial_{t}r(x,t) = e^{tA}Ar(x,0) - i\left[q_1(x) - q_2(x)\right]v(x,t).
\end{split}
\end{equation}
From \eqref{Partial_t_Second_Order_Difference}, note that
\begin{multline} \label{Laplacian_initial_condition}
\Delta r(x,0) =i\sum_{m = 0}^{k}\binom{k}{m}\cdot
\partial^{m}_{z_1}\partial^{k - m}_{z_2}N(0,0) \cdot \\\Big\{\Delta \left[q_{1}(x) - q_{2}(x)\right] \cdot \left[ f(x) \right]^{k} + 2k\left[f(x)\right]^{k-1} \nabla f(x) \cdot \nabla\left[q_{1}(x) - q_{2}(x)\right] \\
+\; k(k-1)\left[f(x)\right]^{k-2}\cdot\left|\nabla f(x) \right|^{2}\cdot\left[q_{1}(x) - q_{2}(x)\right] + k\left[f(x)\right]^{k-1}\cdot \Delta f(x) \cdot [q_{1}(x) - q_{2}(x)]\Big\}.
\end{multline}
Using the estimates in \eqref{linear_well-posedness_result_estimates}, the estimates in \eqref{Partial_t_Second_Order_Difference_source_term_L2_estimate}, equation \eqref{Laplacian_initial_condition}, the Banach algebra property of $H^{2}(\Omega)$, and the Sobolev Embedding Theorem, we obtain the estimates
\begin{equation}
\begin{split}
\Vert r(\cdot\,,t) \Vert_{L^{2}(\Omega)} \leq \Vert r(\cdot\,,0) \Vert_{L^{2}(\Omega)} + T \cdot \sup_{t \in (-T, T)}\Vert \left[(q_1 - q_2)v\right] (\cdot\,,t) \Vert_{L^{2}(\Omega)}  \\
\leq C(\Omega,M,T)\Vert q_{1} - q_{2} \Vert_{L^{2}(\Omega)}\Vert f\Vert^{k}_{H^{4}(\Omega)},\\
\Vert r(\cdot\,,t) \Vert_{H^{2}(\Omega)} \leq C(\Omega, M)\left(\Vert r(\cdot\,,0) \Vert_{H^{2}(\Omega)} + T \cdot \sup_{t \in (-T, T)}\Vert \left[(q_1 - q_2)v\right](\cdot\,,t) \Vert_{H^{2}(\Omega)}\right)\\ \leq
C(\Omega, M, T)\Vert q_{1} - q_{2} \Vert_{H^{2}(\Omega)}\Vert f \Vert^{k}_{H^{4}(\Omega)}
\end{split}
\end{equation}
and
\begin{equation}
\begin{split}
\Vert \partial_{t}r(\cdot\,,t) \Vert_{L^{2}(\Omega)} \leq \Vert (\Delta + p(x))r(\cdot\,,0) \Vert_{L^{2}(\Omega)} +\Vert \left[(q_1 - q_2)v\right](\cdot\,,t) \Vert_{L^{2}(\Omega)} \\ \leq
C(\Omega, M)\Vert q_{1} - q_{2} \Vert_{H^{2}(\Omega)}\Vert f \Vert^{k}_{H^{4}(\Omega)},\\
\Vert \partial_{t}r(\cdot\,,t) \Vert_{H^{2}(\Omega)} \leq C(\Omega, M)\Vert (\Delta + p(x))r(\cdot\,,0) \Vert_{H^{2}(\Omega)} + \Vert \left[(q_1 - q_2)v\right](\cdot\,,t) \Vert_{H^{2}(\Omega)} \\ \leq
C(\Omega, M)\Vert q_{1} - q_{2} \Vert_{H^{4}(\Omega)}\Vert f \Vert^{k}_{H^{4}}.
\end{split}
\end{equation}
\end{proof}
\begin{remark} \label{symmetric_extension_remark}
 We note that under the assumptions of Theorem \ref{main_theorem}-\ref{main_theorem_partial_data}, one may express ${^{(1)}}u_{\epsilon}(x,t) = \overline{{^{(1)}}u_{\epsilon}(x,-t)}$ and $r(x, t) = -\overline{r(x,-t)}$ for $t \in (-T, 0)$, which is justified in Appendix \ref{Extension_Discussion} for completeness. As a result, the estimates in \eqref{L2_H2_Partial_t_Second_Order_Difference_estimates_statement} also hold for $t \in (-T, 0)$. Moreover, for $S \subseteq \partial \Omega$, this extension allows us to control the $L^{2}(S \times (-T, T))$ norm of the Neumann data in terms of its $L^{2}(S \times (0, T))$ norm. This will be used in equation \eqref{initial2_stability_est} concerning the proof of Theorem \ref{main_theorem}.
\end{remark}
We can now proceed to prove Theorem \ref{main_theorem}-\ref{main_theorem_partial_data}. We begin by assuming that $p_{1} \equiv p_{2} \equiv p$ on $\Omega$ and that $q_{1} \equiv q_{2}$ on $\omega$ in order to address \eqref{main_theorem_stable_determination} and \eqref{main_theorem_partial_stable_determination}. When applying the Carleman estimates from Section \ref{main_previous_literature_Carleman_estimate} in the following calculations, the reader should bear in mind that the resulting constants depend on $\Gamma_{0}$.
\subsection{Measurement on a Boundary Subset Satisfying a Geometric Condition (Theorem \ref{main_theorem})} \label{main_theorem_proof_geometric_condition}
The proof of Theorem \ref{main_theorem} now proceeds along the lines of Theorem 2 in \cite{StartingCarlemanEstimate} beginning with an estimate for the initial data in \eqref{Partial_t_Second_Order_Difference}. Recalling the definition of $R_{1}(\cdot)$ from Section \ref{main_previous_literature_Carleman_estimate}, we consider
\begin{equation}\label{initial_condition_estimate02}
\begin{split}
I = \text{Im}\int_{-T}^{0}\int_{\Omega}\left[e^{-s\phi(x,t)}R_1\left(r(x,t)\right)\right]  e^{-s\phi(x,t)}\overline{r(x,t)} \,dx\,dt \\ = 
\text{Im}\int_{-T}^{0}\int_{\Omega} \left[\left(i\partial_t + \Delta + s^{2}\left|\nabla \phi(x,t)\right|^{2}\right)\left(e^{-s\phi(x,t)}r(x,t)\right)\right]  e^{-s\phi(x,t)}\overline{r(x,t)} \,dx\,dt \\ = 
\text{Im}\int_{-T}^{0}\int_{\Omega} \left[\left(i\partial_t + \Delta \right)\left(e^{-s\phi(x,t)}r(x,t)\right)\right] \cdot e^{-s\phi(x,t)}\overline{r(x,t)} \,dx\,dt.
\end{split}
\end{equation}
Note that since $r(x,t) = 0$ on $\partial \Omega \times (0,T)$ from \eqref{Partial_t_Second_Order_Difference}, we have upon integrating by parts that
\begin{equation}
\text{Im}\int_{-T}^{0}\int_{\Omega} \Delta \left(e^{-s\phi(x,t)}r(x,t)\right) \cdot e^{-s\phi(x,t)}\overline{r(x,t)} \,dx\,dt = 0.
\end{equation}
Now, using the fact that $\lim_{t\, \rightarrow\; -T^{+}}e^{-s\phi(x,t)} = 0$ one obtains
\begin{align}\label{integral_calculation}
\begin{split}
I = \text{Re}\int_{-T}^{0}\int_{\Omega}\partial_t\left(e^{-s\phi(x,t)}r(x,t)\right)e^{-s\phi(x,t)}\overline{r(x,t)}\,dx\,dt \\ = 
\frac{1}{2}\int_{-T}^{0}\int_{\Omega}\partial_t\left(e^{-2s\phi(x,t)}|r(x,t)|^2\right) \,dx\,dt
=\frac{1}{2}\int_{\Omega} e^{-2s\phi(x,0)}|r(x,0)|^2\,dx
\\ =
\frac{1}{2}\int_{\Omega} e^{-2s\phi(x,0)} \left|\sum_{m = 0}^{k}\binom{k}{m} \cdot 
\partial^{m}_{z_1}\partial^{k - m}_{z_2}N(0,0)\right|^{2} \cdot \left|q_1(x) - q_2(x)\right|^{2} \cdot \left|f(x)\right|^{2k} \,dx.
\end{split}
\end{align}\\
If we now apply the Cauchy-Schwarz Inequality, the Carleman estimate in \eqref{OldCarlemanEstimate_1} (with $T_{1} = T$), and equation \eqref{boundary_term_bound} to \eqref{initial_condition_estimate02}, we obtain the estimate
\begin{multline} \label{expedited_almost_stability}
I = \frac{1}{2}\int_{\Omega} e^{-2s\phi(x,0)} \left|\sum_{m = 0}^{k}\binom{k}{m} \cdot 
\partial^{m}_{z_1}\partial^{k - m}_{z_2}N(0,0)\right|^{2} \cdot \left|q_1(x) - q_2(x)\right|^{2} \cdot \left|f(x)\right|^{2k} \,dx \\
\leq \left(\int_{-T}^{T}\int_{\Omega}e^{-2s\phi(x,t)}|R_1(r(x,t))|^2 \,dx\,dt\right)^{\frac{1}{2}}\left(\int_{-T}^{T}\int_{\Omega}e^{-2s\phi(x,t)}|r(x,t)|^2\right)^{\frac{1}{2}}\\
\leq \Bigg(C \Vert e^{-s\phi}(q_1 - q_2)v \Vert^{2}_{L^2(\Omega \times (-T,T))}  +
Cs\left\Vert\frac{\partial r}{\partial \nu}\right\Vert^2_{L^2(\Gamma_{0} \times (-T,T))}\Bigg)^{\frac{1}{2}}\cdot \\\Bigg(\frac{C}{s^3}\Vert e^{-s\phi}(q_1 - q_2)v\Vert^{2}_{L^2(\Omega \times (-T,T))} + \frac{C}{s^2}\left\Vert\frac{\partial r}{\partial \nu}\right\Vert^{2}_{L^2(\Gamma_{0} \times (-T,T))}\Bigg)^{\frac{1}{2}} \\ = 
\frac{1}{s^{3/2}}\Bigg(C \Vert e^{-s\phi}(q_1 - q_2)v \Vert^{2}_{L^2(\Omega \times (-T,T))} + Cs\left\Vert\frac{\partial r}{\partial \nu}\right\Vert^{2}_{L^2(\Gamma_{0} \times (-T,T))}\Bigg),
\end{multline}
which holds for $s \in [s_{0}, \infty)$, $\lambda \in [\Lambda_{0}, \infty)$ (where we take $\Lambda_{0} > 1$), and some constant $C = C(\Omega, T, M, s_{0}, \lambda, \Lambda_{0}) > 0$. Let us take $s_{0}$ sufficiently large to satisfy $(1/s)^{3/2} < 1/s$ for $s \in [s_{0},\infty)$. Now, applying the Sobolev Embedding Theorem, the estimates in Lemma \ref{L2_H2_Partial_t_Second_Order_Difference_estimates}, and the fact that $e^{-s\phi(x,t)} \leq e^{-s\phi(x,0)}$ yields the estimate
\begin{equation} \label{Carleman_initial_Condition_Est}
\begin{split}
\int_{\Omega} e^{-2s\phi(x,0)} \left|q_1(x) - q_2(x)\right|^{2} \left|f(x)\right|^{2k} \,dx   \\ \leq \frac{C}{s} \Vert e^{-s\phi(\cdot, 0)}(q_1 - q_2) \Vert^{2}_{L^2(\Omega)} \Vert v \Vert^{2}_{L^{\infty}(-T, T; H^{2}(\Omega))} + 
C\left\Vert\frac{\partial r}{\partial \nu}\right\Vert^{2}_{L^2(\Gamma_{0} \times (-T,T))}\\
\leq \frac{C}{s} \Vert e^{-s\phi(\cdot, 0)}(q_1 - q_2) \Vert^{2}_{L^2(\Omega)} \Vert f \Vert^{2k}_{H^{4}(\Omega)} + 
C\left\Vert\frac{\partial r}{\partial \nu}\right\Vert^{2}_{L^2(\Gamma_{0} \times (-T,T))} \\\leq \frac{C_{1}}{s} \Vert e^{-s\phi(\cdot, 0)}(q_1 - q_2) \Vert^{2}_{L^2(\Omega)} + 
C\left\Vert\frac{\partial r}{\partial \nu}\right\Vert^{2}_{L^2(\Gamma_{0} \times (-T,T))}
\end{split}
\end{equation}
for $s \in [s_{0}, \infty)$, $\lambda \in [\Lambda_{0}, \infty)$, and some constant $C_{1} = C_{1}(\Omega, T, M, \gamma_{+}, s_{0}, \lambda, \Lambda_{0}) > 0$. Using the assumption in Theorem \ref{main_theorem} that $|f(x)| \geq \gamma_{-} > 0$ on $\Omega\setminus\omega$, we have for some constant $C_{2} = C_{2}(\gamma_{-}) > 0$ that the LHS of \eqref{Carleman_initial_Condition_Est} is bounded below as follows:
\begin{equation} \label{LHS_lower_bound_Carleman_Initial_Condition_Est}
\int_{\Omega} e^{-2s\phi(x,0)} \left|q_1(x) - q_2(x)\right|^{2}\left|f(x)\right|^{2k} \,dx \geq C_{2}\int_{\Omega \setminus \omega} e^{-2s\phi(x,0)}\left|q_1(x) - q_2(x)\right|^2\,dx.
\end{equation}
We see that by taking $s = s(C_{1}, C_{2}) \in [s_{0}, \infty)$ sufficiently large, we can absorb the first term on the RHS of \eqref{Carleman_initial_Condition_Est} into the LHS to obtain
\begin{equation} \label{initial2_stability_est_0}
\int_{\Omega} e^{-2s\phi(x,0)} \left|q_1(x) - q_2(x)\right|^{2} \,dx \leq C \left\Vert\frac{\partial r}{\partial \nu}\right\Vert^{2}_{L^2(\Gamma_{0} \times (-T,T))}
\end{equation}
for some $C = C(\Omega, T, M, \gamma_{-}, \gamma_{+}, s_{0}, \lambda, \Lambda_{0}) > 0$. Moreover, note that $\phi(x,0) \leq e^{2\lambda \Vert \psi \Vert_{\infty}}/T^{2}$ so that $e^{-2s\phi(x,0)}$ is bounded below by some constant $C = C(\Omega, T,\lambda, s) > 0$ for $s\in [s_{0},\infty)$.
Using this observation, \eqref{LHS_lower_bound_Carleman_Initial_Condition_Est}, and Remark \ref{symmetric_extension_remark}, equation \eqref{initial2_stability_est_0} becomes
\begin{equation} \label{initial2_stability_est}
\int_{\Omega \setminus \omega} |q_1(x) - q_2(x)|^2\,dx \leq C\int_{\Gamma_{0} \times (-T,T)} \left|\frac{\partial r}{\partial \nu}\right|^2\,dS\,dt \leq C\int_{\Gamma_{0} \times (0,T)} \left|\frac{\partial r}{\partial \nu}\right|^2\,dS\,dt
\end{equation}
for another constant $C = C(\Omega, M, T, \gamma_{-}, \gamma_{+}) > 0$, where we took $\lambda = \Lambda_{0}$ and omitted the parameters $s_{0}, \Lambda_{0}$ due their dependence on the parameters $\Omega, T$ (see Proposition \ref{standard_Schrodinger_Carleman_estimate_1} and Remark \ref{Carleman_weight_observation_at_zero_technical_point}) as well as the parameter $s$ due to its dependence on the constants $C_{1} = C_{1}(\Omega, T, M, \gamma_{+}, s_{0}, \lambda, \Lambda_{0}), C_{2} = C_{2}(\gamma_{-})$ appearing in \eqref{Carleman_initial_Condition_Est}-\eqref{LHS_lower_bound_Carleman_Initial_Condition_Est}. In fact, we have that
\begin{equation} \label{derivative_interchange}
\frac{\partial r}{\partial \nu} = \partial_{\nu}\partial_{t}\left({^{(k)}}u_{\epsilon,1}-{^{(k)}}u_{\epsilon,2}\right) = \partial_t\left(\partial^{k}_{\epsilon}\partial_{\nu}(u_{\epsilon,1}-u_{\epsilon,2})\Big|_{\epsilon = 0}\right).
\end{equation}
This follows by considering the Banach space-valued mapping
\begin{align*}
t \longmapsto \left({^{(k)}}u_{\epsilon,1}-{^{(k)}}u_{\epsilon,2}\right)(\cdot\,,t)\Big|_{\partial \Omega} \in H^{3/2}(\partial \Omega)
\end{align*}
in the weak sense and passing derivatives onto test functions from $C^{\infty}(\partial \Omega)$. This completes the proof of \eqref{main_theorem_stable_determination} in Theorem \ref{main_theorem}.
\begin{remark} \label{derivative_interchange_and_main_theorem_statement}
Observe from the equality in \eqref{derivative_interchange} that the condition $(\partial u_{\epsilon,1}/\partial \nu)\vert_{\Gamma_{0}} = (\partial u_{\epsilon,2}/\partial \nu)\vert_{\Gamma_{0}}$ for all $\epsilon \in (0, \epsilon_{0})$ for some $\epsilon_{0} > 0$ (as in part (c) of Theorem \ref{main_theorem}) implies that $\partial r/\partial \nu = 0$.
\end{remark}
\begin{remark} \label{absorbing_terms_initial_condition_remark}
Since $f(x) \in H^{4}(\Omega)\cap H^{3}_{0}(\Omega)$, $|f(x)|$ becomes arbitrarily small near $\partial \Omega$. We thus need a positive lower bound for $|f(x)|$ away from $\partial \Omega$ in order to absorb terms in \eqref{Carleman_initial_Condition_Est}.
\end{remark}
\textit{Recovery of $p_{j}(x)\;(j = 1,2)$.}\\
The proof of \eqref{main_theorem_stable_determination_0} is similar and has already been treated in \cite{StartingCarlemanEstimate}. Observe from \eqref{FirstLinear_for_p1} that 
\begin{equation} \label{First_Linear_Difference}
\begin{cases} \left(i\partial_t + \Delta + p_1(x)\right)\left({^{(1)}}u_{\epsilon,1} - {^{(1)}}u_{\epsilon,2}\right)(x,t) = \left[p_2(x) - p_1(x)\right]{^{(1)}}u_{\epsilon,2}(x,t) \text{ on } \Omega \times (0,T),\\ \left({^{(1)}}u_{\epsilon,1} - {^{(1)}}u_{\epsilon,2}\right)(x,t) = 0 \text{ on } \partial \Omega \times (0,T),\\ \left({^{(1)}}u_{\epsilon,1} - {^{(1)}}u_{\epsilon,2}\right)(x,0) = 0 \text{ on } \Omega \times \{0\}.\end{cases}
\end{equation}
Moreover, proceeding as in Lemma \ref{Partial_t_Second_Order_Difference_lemma}, we can show that $r_{1}(x,t) = \partial_{t}\left({^{(1)}}u_{\epsilon,1} - {^{(1)}}u_{\epsilon,2}\right)(x,t)$ satisfies
\begin{equation} \label{Partial_t_First_Linear_Difference}
\begin{cases} \left(i\partial_t + \Delta + p_1(x)\right)r_{1}(x,t) = \left[p_2(x) - p_1(x)\right]\partial_{t}{^{(1)}}u_{\epsilon,2}(x,t) \text{ on } \Omega \times (0,T),\\ r_{1}(x,t) = 0 \text{ on } \partial \Omega \times (0,T),\\ 
r_{1}(x,0) = -i\left[p_2(x) - p_1(x)\right]f(x) \text{ on } \Omega \times \{0\}.\end{cases}
\end{equation}
Now, repeating the steps used to prove \eqref{main_theorem_stable_determination} yields \eqref{main_theorem_stable_determination_0}. In view of Remark \ref{derivative_interchange_and_main_theorem_statement}, combining \eqref{main_theorem_stable_determination_0} and \eqref{main_theorem_stable_determination} yields unique determination of $p_{j}(x),q_j(x)\;(j = 1,2)$ assuming knowledge of the coefficients near $\partial \Omega$.

\subsection{Measurement on an Arbitrary Subset of the Boundary (Theorem \ref{main_theorem_partial_data})}
In order to realize stable determination from a measurement on an arbitrary relatively open subset $\Gamma \subset \partial \Omega$, the authors of \cite{MainMethod} prove a unique continuation estimate\footnote{The authors of \cite{MainMethod} refer to this estimate as an ``observability inequality".} that relates the $H^{1}(\Omega_{2} \setminus \Omega_{1})$-norm in the Carleman estimate from Proposition \ref{standard_Schrodinger_Carleman_estimate_0} to the Neumann data on $\Gamma$. That is, keeping in mind the assumptions in Theorem \ref{main_theorem_partial_data}, we have
\begin{proposition} \label{logarithmic_observability_inequality_proposition}
There exist constants $\gamma_{*}, C, \mu > 0$ and $T_{1} \in (0, T)$ such that the solution $r(x,t)$ to \eqref{Partial_t_Second_Order_Difference} satisfies
\begin{equation} \label{logarithmic_observability_inequality_statement}
\int_{\Omega_{2}\setminus \Omega_{1}}\int_{-T_{1}}^{T_{1}}\left|\nabla r(x,t) \right|^{2} + \left|r(x,t) \right|^{2}\,dt\,dx \leq C\left(\frac{1}{\gamma} + e^{-\mu \gamma} +  e^{\mu \gamma}\int_{\Gamma}\int_{-T}^{T}  \left|\frac{\partial r}{\partial \nu}(x,t)\right|^{2} \,dt \,dS \right)
\end{equation}
for $\gamma \geq \gamma_{*}$. Recall that $\Omega_{j}\; (j = 1,2)$ is defined by \eqref{partial_data_proof_subset_1}.
\end{proposition}
In turn, the proof of the unique continuation estimate in \eqref{logarithmic_observability_inequality_statement} relies on a parabolic Carleman estimate that enables an estimate by the Neumann data on $\Gamma$. Recalling that $p_{1} \equiv p_{2} \equiv p$ on $\Omega$ and that $q_{1} \equiv q_{2}$ on $\omega$, we shall follow the argument in \cite{MainMethod} in relation to \eqref{Partial_t_Second_Order_Difference}. \\ \indent
For a fixed parameter $h \in (0,1)$, the authors of \cite{MainMethod} establish a Carleman estimate for the parabolic operator $h^{-1}\partial_{\tau} - \Delta - p(x)$. Before stating it, we need to introduce a function $\psi_{0}(x)$ defined as follows.
\begin{proposition} \label{parabolic_weight_function_0}
Let $\Gamma$ be a relatively open subset of $\partial \Omega$ and consider $\omega \subset \Omega$ a relatively open neighborhood of $\partial \Omega$ such that $\partial\omega \setminus \partial\Omega$ is $C^{2}$. Then, there exists $\psi_{0}(x) \in C^{2}(\overline{\omega})$ satisfying
\begin{align}
\begin{cases}
\psi_{0}(x) > 0 \quad (x \in \omega),\\
\left| \nabla \psi_{0}(x)\right| > 0 \quad (x \in \overline{\omega}),\\
\psi_{0}(x) = 0 \quad (x \in \partial \omega \setminus \Gamma),\\
\frac{\partial \psi_{0}}{\partial \nu}(x) \leq 0 \quad (x \in \partial \omega \setminus \Gamma).
\end{cases}
\end{align}
\end{proposition}
\begin{remark} \label{smoothness_Pi_inner_boundary}
The existence of the function defined in the preceding proposition is proved through Lemma $2.1$ and Lemma $2.3$ in \cite{IY1998}. The proof given in \cite{IY1998} uses the regularity assumption on $\partial\omega \setminus \partial\Omega$ to ensure that we can write $\nabla \psi_{0}(x) = (\partial \psi_{0}/\partial\nu)(x)\nu(x)$ on $\partial \Omega \setminus \Gamma$ after we know that $\psi(x) = 0$ on $\partial \Omega \setminus \Gamma$.
\end{remark}
Now, we recall below (with minor changes in notation) the parabolic Carleman estimate from Lemma $4.2$ in \cite{MainMethod}.
\begin{proposition} \label{parabolic_Carleman_estimate_0}
Let $p(x) \in L^{\infty}(\Omega)$ with $\Vert p\Vert_{L^{\infty}(\Omega)} \leq M$ for some $M > 0$.  Letting $\psi_{0}(x)$ be defined as in Proposition \ref{parabolic_weight_function_0} and denoting $\Vert \psi_{0}\Vert_{\infty}$ = $\Vert \psi_{0}\Vert_{L^{\infty}(\Omega)}$, we define
\begin{equation} \label{auxilliary_weight}
\theta_{0}(x,\tau) = \frac{e^{\lambda \psi_{0}(x)}}{1 - \tau^{2}} \quad (x \in \omega,\; \tau \in (-1,1)).
\end{equation}
Picking $\Vert \psi_{0}\Vert_{L^{\infty}(\Omega)} < a < b < 2a - \Vert \psi_{0}\Vert_{L^{\infty}(\Omega)}$, we also define
\begin{equation} \label{weight_function_exponent_0}
\phi_{0}(x,\tau) = \frac{e^{\lambda(\Vert \psi_{0}\Vert_{\infty} + b)} - e^{\lambda(\psi_{0}(x) + a)}}{1 - \tau^{2}} \quad (x \in \omega,\; \tau \in (-1,1)).
\end{equation}
For $h \in (0, 1)$ a fixed parameter, there exist $\lambda_{*} = \lambda_{*}(\omega) > 0$, $\sigma_{*} = \sigma_{*}(\omega) > 0$ such that for all $\lambda \geq \lambda_{*}$, $\sigma \geq \sigma_{*}/h$, and a constant $C = C(\omega, \Gamma, \sigma_{*},\lambda, \lambda_{*}, h, M) > 0$, we have
\begin{equation} \label{parabolic_Carleman_estimate_statement}
\begin{split}
\sigma^3 \Vert e^{-\sigma\phi_{0}}w\Vert_{L^2(\omega \times (-1,1))}^2 + \sigma\Vert e^{-\sigma\phi_{0}}\nabla w\Vert_{L^2(\omega \times (-1,1))}^2 \\
\leq C\left(\Vert e^{-\sigma\phi_{0}} (h^{-1}\partial_{\tau} - \Delta - p(x))w\Vert_{L^2(\omega \times (-1,1))}^2 + \sigma\left\Vert \sqrt{\theta_{0}}e^{-\sigma\phi_{0}}\frac{\partial w}{\partial \nu} \right\Vert_{L^2(\Gamma \times (-1,1))}^2\right).
\end{split}
\end{equation}
The above holds for $w \in L^{2}\left(-1,1;\,H^{1}_{0}(\omega) \right)$ such that $(h^{-1}\partial_{\tau} - \Delta - p(x))w \in L^{2}\left(\omega \times (-1,1)\right)$ and $\partial w/\partial \nu  \in L^{2}\left(-1,1;\,L^{2}(\Gamma) \right)$.
\end{proposition}
\begin{remark} \label{partial_data_Carleman_estimate_parameter_dependence}
The dependence on $\omega$ and $\Gamma$ for the constant in \eqref{parabolic_Carleman_estimate_statement} is not explicitly stated in Lemma $4.2$ from \cite{MainMethod}. However, this dependence is apparent from the proof provided in Appendix A of \cite{MainMethod}, where the constant eventually depends on bounds for the function $\theta_{0}$ (defined in \eqref{auxilliary_weight}) over $\omega \times (-1, 1)$. Moreover, $\theta_{0}$ is defined in terms of the function $\psi_{0}$ given in Proposition \ref{parabolic_weight_function_0}, which yields the dependence on $\Gamma$.
\end{remark}
In a series of steps, we proceed as in \cite{MainMethod} to relate the solution to the Schr\"odinger IBVP \eqref{Partial_t_Second_Order_Difference} to the solution of a certain parabolic equation in order to make use of the partial boundary measurement appearing in \eqref{parabolic_Carleman_estimate_statement}, consequently proving Proposition \ref{logarithmic_observability_inequality_proposition}. For the following calculations, whenever we apply the Carleman estimate from Proposition \ref{parabolic_Carleman_estimate_0}, it should be understood that the resulting constants depend on $\Gamma$.\\ \\
\textit{Step 1 (Transforming \eqref{Partial_t_Second_Order_Difference} into a parabolic IBVP on $\omega$).}\\ 
In order to apply the Carleman estimate \eqref{parabolic_Carleman_estimate_statement}, we need to transform \eqref{Partial_t_Second_Order_Difference} into a parabolic IBVP on $\omega \times (-1,1)$ with vanishing boundary conditions on $\partial \omega \times (-1,1)$. At the same time, it turns out to be convenient to have a positive lower bound for the parabolic Carleman weight, namely $\phi_{0}(x,\tau)$ defined in \eqref{weight_function_exponent_0}, in order to eventually obtain a stability estimate.\\ \indent 
We begin by obtaining the analogues of equations $(4.13)$ and $(4.20)-(4.22)$ in \cite{MainMethod}. First, note that we can express $\partial \omega = \Pi \cup \partial \Omega$. Using Proposition \ref{parabolic_weight_function_0}, we can find a constant $\kappa > 0$ such that
\begin{equation} \label{weight_function_lower_bound_0}
\psi_{0}(x) \geq 2\kappa \quad (x \in \Omega_{2}\setminus \Omega_{1} = \overline{\omega_{1}} \setminus \overline{\omega_{2}}).
\end{equation}
Moreover, since $\psi_{0}(x) = 0$ on $\Pi$, we can find an open neighborhood $\Pi' \subset \omega$ of $\Pi$ such that 
\begin{equation}
\overline{\Pi'} \cap \overline{\omega_{1}} = \emptyset
\end{equation}
and 
\begin{equation} \label{weight_function_upper_bound_0}
\psi_{0}(x) \leq \kappa \quad (x \in \Pi').
\end{equation}
Now, taking an open neighborhood $\Pi'' \subset \Pi'$ of $\Pi$ (such that $\overline{\Pi''} \cap (\overline{\omega \setminus \Pi'}) = \emptyset$), define a smooth function $\chi(x)$ for $x \in \omega$ such that
\begin{equation} \label{smooth_cutoff_0}
\chi(x) = \begin{cases}0,\; x \in \Pi''\\
1,\; x \in \omega \setminus \Pi'.
\end{cases}
\end{equation}
Recalling that $r(x,t)$ satisfies \eqref{Partial_t_Second_Order_Difference} with source term $[q_{1}(x)-q_{2}(x)]v(x,t)$ (where $v(x,t)$ is defined by \eqref{Partial_t_Second_Order_Difference_source_term}), we have using the fact that $q_{1}(x) = q_{2}(x)\; (x \in \omega)$ and $r(x,0) = 0\; (x \in \omega)$ that $\tilde{r}(x,t) = \chi(x)r(x,t)$ satisfies
\begin{equation} \label{Partial_t_Second_Order_Difference_near_boundary}
\begin{cases}
(i\partial_{t} + \Delta + p(x))\tilde{r}(x,t) = r(x,t)\Delta\chi(x) + 2\left(\nabla \chi(x) \cdot \nabla r(x,t) \right)\\
= \left[\Delta, \chi(x) \right]r(x,t) \text{ on } \omega \times (-T,T),\\
\tilde{r}(x,t) = 0 \text{ on } \partial \omega \times (-T, T),\\
\tilde{r}(x,0) = 0 \text{ on } \omega \times \{0\}.
\end{cases}
\end{equation}
We have used Remark \ref{symmetric_extension_remark} to consider the above equation on $(-T, T)$. Observe from \eqref{smooth_cutoff_0} that the source term $\left[\Delta, \chi(x) \right]r(x,t)$ in \eqref{Partial_t_Second_Order_Difference_near_boundary} is supported on $\Pi'$. \\ \indent
Now, we transform \eqref{Partial_t_Second_Order_Difference_near_boundary} into a parabolic IBVP using a variant of the Fourier-Bros-Iagolnitzer (F.B.I.) transformation $\mathcal{F}_{\gamma}$, defined as in \cite{MainMethod} by
\begin{equation} \label{Fourier_Bros_Iagnolnitzer_transform_0}
\begin{split}
w_{\gamma}(x,z) = w_{\gamma, t}(x,\tau) = 
\mathcal{F}_{\gamma}(w(x,t)) = \sqrt{\frac{\gamma}{2\pi}} \int_{\mathbb{R}}e^{-(\gamma/2)(z - \eta)^{2}}\theta(\eta)w(x,\eta h)\,d\eta,
\end{split}
\end{equation}
where $z = t - i\tau$. In the above, $T_{0} > T > 0$ and $\gamma > 0$ are parameters that will be fixed later, $\tau \in (-1,1)$,
\begin{equation} \label{time_smooth_cutoff}
\theta(\eta) = \begin{cases}
0,\; |\eta| \geq 3T_{0} \\
1,\; |\eta| \leq 2T_{0}
\end{cases}
\end{equation}
is a smooth cut-off function, and $h = T/3T_{0} \in (0,1)$. For $w(x,t) \in L^{2}(-T, T;H^{1}(\omega))$, we have $w_{\gamma, t}(x,\tau) \in L^{2}(-1, 1; H^{1}(\omega))$. Likewise, for $w(x,t) \in L^{2}(\omega \times (-T,T))$, we have $w_{\gamma, t}(x,\tau) \in L^{2}(\omega \times (-1,1))$. We also note that $w_{\gamma}(x,z)$ is analytic with respect to $z = t - i\tau$.\\ \indent
Using \eqref{Partial_t_Second_Order_Difference_near_boundary}, one can show that $r_{0}(x, t-i\tau) = r_{0}(x,z) = \mathcal{F}_{\gamma}(\tilde{r}(x,t))$ satisfies
\begin{equation} \label{Partial_t_Second_Order_Difference_parabolic}
\begin{cases}
(h^{-1}\partial_{\tau} - \Delta - p(x))r_{0}(x,z) = -\sqrt{\frac{\gamma
}{2\pi}}\int_{\mathbb{R}}e^{-(\gamma/2)(z - \eta)^{2}}\theta(\eta)\left[\Delta, \chi(x)\right]r(x,\eta h)\,d\eta\\
-\frac{i}{h}\sqrt{\frac{\gamma}{2\pi}}\int_{\mathbb{R}}e^{-(\gamma/2)(z - \eta)^{2}}\theta'(\eta)r(x,\eta h)\,d\eta \text{ on } \omega \times (-1,1),\\
r_{0}(x,z) = 0 \text{ on } \partial\omega \times (-1,1),\\
r_{0}(x,z) = r_{0}(x,t) \text{ on } \omega \times \{0\}.
\end{cases}
\end{equation}
Applying the Carleman estimate \eqref{parabolic_Carleman_estimate_statement} to $r_{0}(x,z)$ satisfying \eqref{Partial_t_Second_Order_Difference_parabolic}, one can proceed as in the proofs of Lemma $4.3$ and $4.4$ in \cite{MainMethod} to show that there exists $\epsilon = \epsilon(\omega) \in (0, 1)$ such that
\begin{equation} \label{stability_estimate_step_3_0_0_0}
\begin{split}
\int_{-\epsilon}^{\epsilon}\int_{\Omega_{2}\setminus \Omega_{1}}\left( |\nabla r_{0}(x,z)|^{2} + |r_{0}(x,z)|^{2}\right)\,dx\,d\tau \leq \\
\frac{C e^{2\gamma}e^{-\mu_{0}\sigma}}{\sigma}\int_{\omega}\int_{-3T_{0}}^{3T_{0}} \left|\left[\Delta, \chi(x)\right]r(x,\eta h) \right|^{2} \,d\eta\,dx \\
+\; \frac{C e^{2\gamma-\gamma T_{0}^{2}}e^{\tilde{\mu}\sigma}}{h^{2}\sigma}\int_{\omega}\int_{-3T_{0}}^{3T_{0}}\left|r(x,\eta h) \right|^{2}\,d\eta\,dx \\
+\;  e^{2\gamma}e^{\tilde{\mu}\sigma}C(T_{0})\int_{-1}^{1}\int_{\Gamma}\int_{-3T_{0}}^{3T_{0}} \theta_{0}(x,\tau)e^{-2\sigma\phi_{0}(x,\tau)} \left|\frac{\partial r}{\partial \nu}(x,\eta h)\right|^{2} \,d\eta \,dS\,d\tau
\end{split}
\end{equation}
for $\sigma \geq \sigma_{*}/h > 1$, some constant $\mu_{0} > 0$, and
\begin{equation} \label{parabolic_Carleman_estimate_boundary_term_constant_0}
\tilde{\mu} = 2\frac{e^{\lambda(\Vert \psi_{0}\Vert_{\infty} + b)} - e^{\lambda(2\kappa + a)}}{1 - \epsilon^{2}}.
\end{equation}
In writing \eqref{stability_estimate_step_3_0_0_0}, we used the fact that
\begin{equation} \label{integrand_estimate_0}
\left|e^{-(\gamma/2)(z-\eta)^{2}} \right|^{2} = e^{-\gamma(t-\eta)^{2}}e^{\gamma \tau^{2}}
\end{equation}
and that $\gamma e^{\gamma} < e^{2 \gamma}$ for $\gamma > 0$.\\ \indent
Now, if we take sufficiently large $T_{0} = T_{0}(\omega) > T > 0$ and consider $\sigma = T_{0}\gamma$, we can show that
\begin{equation} \label{stability_estimate_step_3_0_0}
\int_{-\epsilon}^{\epsilon}\int_{\Omega_{2}\setminus \Omega_{1}}\left( |\nabla r_{0}(x,z)|^{2} + |r_{0}(x,z)|^{2}\right)\,dx\,d\tau \leq 
C\left(e^{-\mu \gamma} +  e^{\mu \gamma}\int_{-1}^{1}\int_{\Gamma} \left| \frac{\partial r_{0}}{\partial \nu}(x,z)\right|^{2}\,dS\,d\tau\right)
\end{equation}
for constants $C = C(\Omega, \omega, \Gamma, M, \tilde{M}, \gamma_{+}, \sigma_{*}, T, T_{0}) > 0$, $\gamma \geq \gamma_{*}(\omega, T_{0}) = \sigma_{*}/T_{0}$, and $\mu = \mu(\omega, T, T_{0}) > 0$, where we have $\sigma_{*} = \sigma_{*}(\omega)$, $T_{0} = T_{0}(\omega) > T$, and $\epsilon = \epsilon(\omega) \in (0,1)$. Using \eqref{integrand_estimate_0} for the boundary integral in \eqref{stability_estimate_step_3_0_0}, we further have
\begin{equation} \label{stability_estimate_step_3_0}
\begin{split}
\int_{-1}^{1}\int_{\Gamma} \left| \frac{\partial r_{0}}{\partial \nu}(x,z)\right|^{2}\,dS\,d\tau \leq 
C\gamma \int_{-1}^{1}\int_{\Gamma} \left| \frac{\partial}{\partial \nu}\left( \int_{-3T_{0}}^{3T_{0}} e^{-(\gamma/2)(z-\eta)^{2}} \theta(\eta) r(x,\eta h) \,d\eta \right)\right|^{2}\,dS\,d\tau \\ \leq
\gamma e^{\gamma}C(T_{0})\int_{-1}^{1}\int_{\Gamma}\int_{-3T_{0}}^{3T_{0}}  \left|\frac{\partial r}{\partial \nu}(x,\eta h)\right|^{2} \,d\eta \,dS\,d\tau.
\end{split}
\end{equation}
Now, if we change variables to integrate with respect to ``$t = \eta h$", we have that
\begin{equation} \label{stability_estimate_step_3}
\int_{-\epsilon}^{\epsilon}\int_{\Omega_{2}\setminus \Omega_{1}}\left( |\nabla r_{0}(x,z)|^{2} + |r_{0}(x,z)|^{2}\right)\,dx\,d\tau \leq 
C\left(e^{-\mu \gamma} +  e^{\mu \gamma}\int_{\Gamma}\int_{-T}^{T}  \left|\frac{\partial r}{\partial \nu}(x,t)\right|^{2} \,dt \,dS\right)
\end{equation}
for $\gamma \geq \gamma_{*}(\omega, T_{0})$ and constants $C = C(\Omega, \omega, \Gamma, M, \tilde{M}, \gamma_{+}, T, T_{0}) > 0$, $\mu = \mu(\omega, T, T_{0}) > 0$. We have omitted ``$\sigma_{*}$" since it depends on $\omega$.
\begin{remark}
We note that $\epsilon \in (0,1)$ is chosen to ensure that the first term on the RHS of \eqref{stability_estimate_step_3} has a negative exponent. It is possible to pick such an $\epsilon$ due to the bounds for the parabolic Carleman weight in \eqref{weight_function_lower_bound_0}, \eqref{weight_function_upper_bound_0}, and the fact that $\left[\Delta, \chi(x) \right]r(x,t)$ in \eqref{Partial_t_Second_Order_Difference_near_boundary} is supported on $\Pi'$.
\end{remark}
\textit{Step 2 (Estimates for $r_{0}(x,z)\big|_{\tau = 0}$ and $\nabla r_{0}(x,z)\big|_{\tau = 0}$).}\\ 
The next step is motivated by the observation that
\begin{equation} \label{Schrodinger_parabolic_relation_0}
r_{0}(x,z)\big|_{\tau = 0} = r_{0}(x,t) = \sqrt{\frac{\gamma}{2\pi}} \int_{\mathbb{R}}e^{-(\gamma/2)(t - \eta)^{2}}\theta(\eta)\tilde{r}(x,\eta h)\,d\eta = (K_{\gamma} * \theta\tilde{r}_{h})(x,t),
\end{equation}
where we have denoted
\begin{equation} \label{Schrodinger_parabolic_relation_1}
K_{\gamma}(\cdot) = \sqrt{\frac{\gamma}{2\pi}}e^{-(\gamma/2)(\cdot)^{2}},\; \theta\tilde{r}_{h}(x,\eta) = \theta(\eta)\tilde{r}(x, h\eta).
\end{equation}
In other words, \eqref{Schrodinger_parabolic_relation_0}-\eqref{Schrodinger_parabolic_relation_1} provide a way of relating the solution to the parabolic IBVP \eqref{Partial_t_Second_Order_Difference_parabolic} to the solution of the Schr\"odinger IBVP \eqref{Partial_t_Second_Order_Difference}. In the following, we work with ``$r_{0}(x,z)$" though similar calculations hold for ``$\nabla r_{0}(x,z)$".\\ \indent
In order to make use of \eqref{Schrodinger_parabolic_relation_0}-\eqref{Schrodinger_parabolic_relation_1}, for each fixed $t \in (-2T/3, 2T/3)$ we apply the Cauchy Integral Formula on the boundary of a disk $B_{\epsilon_{0}}(t-i0)$ with radius $\epsilon_{0} \in (0, \epsilon)$ small enough to ensure
\begin{equation} \label{rectangle_constraint_Cauchy_formula}
z = t -i\tau \in B_{\epsilon_{0}}(t-i0) \implies \text{Re}(z) \in (-T, T), \;\text{Im}(z) \in (-\epsilon, \epsilon).
\end{equation}
We can write
\begin{equation} \label{Cauchy_formula_0}
r_{0}(x, t) = \frac{1}{2\pi i}\int_{\partial B_{\epsilon_{0}}(t-i0)}\frac{r_{0}(x,z)}{z - t}\,dz.
\end{equation}
With respect to the polar coordinates $\rho = |z - t|$ and $\Theta = \text{Arg}(z-t) \in (-\pi, \pi]$, \eqref{Cauchy_formula_0} becomes
\begin{equation} \label{Cauchy_formula_1}
\begin{split}
r_{0}(x,t) = \frac{1}{2\pi}\int_{-\pi}^{\pi} r_{0}(x, t + \rho e^{i \Theta})\,d\Theta \\
\implies r_{0}(x,t) = \frac{1}{2\pi \epsilon_{0}}\int_{0}^{\epsilon_{0}}\int_{-\pi}^{\pi} r_{0}(x, t + \rho e^{i \Theta})\,d\Theta\,d\rho \\
\implies \left|r_{0}(x,t) \right|^{2} \leq C(\omega, T)\int_{-T}^{T}\int_{-\epsilon}^{\epsilon} \left|r_{0}(x,t - i\tau) \right|^{2}\,d\tau\, dt,
\end{split}
\end{equation}
In the last line above, we have majorized the integral over the disk $B_{\epsilon_{0}}(t-i0)$ by the integral over $(-T,T) \times (-\epsilon, \epsilon)$. We also applied the Cauchy-Schwarz Inequality and absorbed the resulting factors into $C(\omega, T)$. Recall that that the choice of $\epsilon_{0}$ depends on $T$ and $\epsilon$, which in turn depends on $\omega$.\\ \indent Now, integrating over $(\Omega_{2}\setminus \Omega_{1}) \times (-2T/3, 2T/3)$ and applying \eqref{stability_estimate_step_3} yield
\begin{equation} \label{stability_estimate_step_4}
\int_{-2T/3}^{2T/3}\int_{\Omega_{2}\setminus\Omega_{1}}\left|r_{0}(x,t) \right|^{2}\,dx\,dt \leq C\left(e^{-\mu \gamma} +  e^{\mu \gamma}\int_{\Gamma}\int_{-T}^{T}  \left|\frac{\partial r}{\partial \nu}(x,t)\right|^{2} \,dt \,dS\right)
\end{equation}
for $\gamma \geq \gamma_{*}(\omega, T_{0})$ and constants $C = C(\Omega, \omega, \Gamma, M, \tilde{M}, \gamma_{+}, T, T_{0}) > 0$, $\mu = \mu(\omega, T, T_{0}) > 0$.
Carrying out the same calculation for ``$\nabla r_{0}(x,z)$" yields
\begin{equation} \label{stability_estimate_4_nabla}
\int_{-2T/3}^{2T/3}\int_{\Omega_{2}\setminus\Omega_{1}}\left|\nabla r_{0}(x,t) \right|^{2}\,dx\,dt \leq C\left(e^{-\mu \gamma} +  e^{\mu \gamma}\int_{\Gamma}\int_{-T}^{T}  \left|\frac{\partial r}{\partial \nu}(x,t)\right|^{2} \,dt \,dS\right).
\end{equation} \\
\textit{Step 3 (Relating ``$r_{0}(x,z)\big|_{\tau = 0}$" and ``$\nabla r_{0}(x,z)\big|_{\tau = 0}$" to \eqref{Partial_t_Second_Order_Difference}).} \\ \indent
Now, we will apply the observation in \eqref{Schrodinger_parabolic_relation_0} to \eqref{stability_estimate_step_4}. In the following, we set $\tilde{r}(\cdot,t)$ equal to zero outside $(-2T/3, 2T/3)$ (and consequently $\tilde{r}_{h}(\cdot, \eta)$ equal to zero outside $(-2T_{0}, 2T_{0})$). Letting $\widehat{g}(x,\zeta)$ denote the Fourier transform of $g(x,\eta)$ with respect to $\eta$ and denoting $(\theta\tilde{r}_{h})(x,\eta) = \theta(\eta)\tilde{r}_{h}(x,\eta) = \theta(\eta)\tilde{r}(x,\eta h)$, we obtain
\begin{equation} \label{Fourier_transform_estimate_0}
\left|\widehat{(\theta\tilde{r}_{h})}(x,\zeta) - \widehat{r_{0}}(x,\zeta h)\right|^{2} = \left|\widehat{(\theta\tilde{r}_{h})}(x,\zeta) - \widehat{K_{\gamma}}(\zeta)\widehat{(\theta\tilde{r}_{h})}(x,\zeta)\right|^{2}
\leq \frac{C\zeta^{2}}{\gamma}\left|\widehat{(\theta\tilde{r}_{h})}(x,\zeta) \right|^{2}
\end{equation}
using the estimate
\begin{equation} \label{Fourier_transform_estimate_1}
\left|1 - \widehat{K_{\gamma}}(\zeta)\right|^{2} \leq \frac{C\zeta^{2}}{\gamma}
\end{equation}
for some constant $C > 0$. We see from Plancherel's Theorem that
\begin{equation}
\begin{split}
\int_{\Omega_{2}\setminus \Omega_{1}}\int_{-2T_{0}}^{2T_{0}} \left|\theta(\eta)\tilde{r}(x,\eta h) - r_{0}(x,\eta h) \right|^{2}\,d\eta\,dx = \int_{\Omega_{2}\setminus \Omega_{1}}\int_{\mathbb{R}} \left|\widehat{(\theta\tilde{r}_{h})}(x,\zeta) - \widehat{r_{0}}(x,\zeta h) \right|^{2}\,d\zeta\,dx\\
\leq \frac{C}{\gamma} \int_{\Omega_{2}\setminus \Omega_{1}}\int_{\mathbb{R}} \left|\zeta\widehat{(\theta\tilde{r}_{h})}(x,\zeta)\right|^{2} \,d\zeta\,dx = \frac{C}{\gamma} \int_{\Omega_{2}\setminus \Omega_{1}}\int_{-2T_{0}}^{2T_{0}} \left|\partial_{\eta}(\theta\tilde{r})(x,h\eta)\right|^{2} \,d\eta\,dx \\
\leq \frac{C(h)}{\gamma} \int_{\Omega_{2}\setminus \Omega_{1}}\int_{-2T/3}^{2T/3} \left|\partial_{t}\left[\theta(t/h)\tilde{r}(x,t)\right]\right|^{2} \,dt\,dx.
\end{split}
\end{equation}
Noting from \eqref{time_smooth_cutoff} that $\partial_{t}\theta(t/h)$ is supported on $[-T, -2T/3] \cup [2T/3, T]$ and using the estimates in Lemma \ref{L2_H2_Partial_t_Second_Order_Difference_estimates}, we bound the last term above by
\begin{equation} \label{Fourier_transform_estimate_2_0}
\frac{C(T,h)}{\gamma}\int_{\Omega_{2}\setminus \Omega_{1}}\int_{-T}^{T} \left|\partial_{t}r(x,t)\right|^{2}  \,dt\,dx \leq \frac{C(\Omega, M, \tilde{M}, \gamma_{+}, T, T_{0})}{\gamma}
\end{equation}
Using \eqref{stability_estimate_step_4}, we have
\begin{equation} \label{Fourier_transform_estimate_2}
\begin{split}
\frac{1}{h}\int_{\Omega_{2}\setminus \Omega_{1}}\int_{-2T/3}^{2T/3}\left|r(x,t) \right|^{2}\,dt\,dx = \int_{\Omega_{2}\setminus \Omega_{1}}\int_{-2T_{0}}^{2T_{0}}\left|\left(\theta\tilde{r}_{h}\right)(x,\eta) \right|^{2}\,d\eta\,dx \leq \\
\int_{\Omega_{2}\setminus \Omega_{1}}\int_{-2T_{0}}^{2T_{0}} \left|\theta(\eta)\tilde{r}(x,\eta h) - r_{0}(x,\eta h) \right|^{2}\,d\eta\,dx + \frac{1}{h}\int_{\Omega_{2}\setminus \Omega_{1}}\int_{-2T/3}^{2T/3} \left|r_{0}(x,t) \right|^{2}\,dt\,dx \\
\implies \int_{\Omega_{2}\setminus \Omega_{1}}\int_{-2T/3}^{2T/3}\left|r(x,t) \right|^{2}\,dt\,dx \leq C\left(\frac{1}{\gamma} + e^{-\mu \gamma} +  e^{\mu \gamma}\int_{\Gamma}\int_{-T}^{T}  \left|\frac{\partial r}{\partial \nu}(x,t)\right|^{2} \,dt \,dS \right)
\end{split}
\end{equation}
for $\gamma \geq \gamma_{*}$ and constants $C = C(\Omega, \omega, \Gamma, M, \tilde{M}, \gamma_{+}, T, T_{0}) > 0$, $\mu = \mu(\omega, T, T_{0}) > 0$. In writing the first equality, we used the fact that $\chi(x) = 1$ for $x \in \Omega_{2}\setminus\Omega_{1}$ and that $\theta(\eta) = 1$ for $\eta \in [-2T_{0}, 2T_{0}]$. Repeating the preceding calculation for ``$\nabla r_{0}(x,t)$" and combining with \eqref{Fourier_transform_estimate_2} yields \eqref{logarithmic_observability_inequality_statement}, which completes the proof of Proposition \ref{logarithmic_observability_inequality_proposition}.
%\begin{equation} \label{stability_estimate_step_5}
%\int_{\Omega_{2}\setminus \Omega_{1}}\int_{-2T/3}^{2T/3}\left|\nabla r(x,t) %\right|^{2} + \left|r(x,t) \right|^{2}\,dt\,dx \leq C\left(\frac{1}{\gamma} + %e^{-\mu \gamma} +  e^{\mu \gamma}\int_{\Gamma}\int_{-T}^{T}  %\left|\frac{\partial r}{\partial \nu}(x,t)\right|^{2} \,dt \,dS \right)
%\end{equation}
%for $\gamma \geq \gamma_{*}$ and constants $C = C(\Omega, \omega, M, \tilde{M}, %\gamma_{+}, T, T_{0}) > 0$, $\mu = \mu(\omega, T, T_{0}) > 0$.
\begin{remark} \label{higher_regularity_partial_data_case}
Recalling Lemma \ref{L2_H2_Partial_t_Second_Order_Difference_estimates}, we observe in relation to \eqref{Fourier_transform_estimate_2_0} that the calculation for ``$r_{0}(x,t)$" only requires a bound on the $H^{2}(\Omega)$-norm of $q_{2}-q_{1}$. The process of obtaining the analogue of \eqref{Fourier_transform_estimate_2} for ``$\nabla r_{0}(x,t)$" leads us to estimate $\Vert \partial_{t}(\nabla r) \Vert_{L^{2}((\Omega_{2}\setminus\Omega_{1}) \times (-T_{1}, T_{1}))}$ in place of \eqref{Fourier_transform_estimate_2_0}. In order to estimate this term with Lemma \ref{L2_H2_Partial_t_Second_Order_Difference_estimates}, we require higher regularity for the coefficients (specifically a bound on the $H^{4}(\Omega)$-norm of $q_{2} - q_{1}$).
\end{remark}
\textit{Completing the proof of Theorem \ref{main_theorem_partial_data}.}\\ 
 Recall $\phi(x,t)$ and the parameter $s > 0$ defined in Proposition \ref{standard_Schrodinger_Carleman_estimate_0}. Denoting $r_{1}(x,t) = e^{-s\phi(x,t)}r(x,t)$ and 
\begin{equation}
\begin{split}
m(x,t) = e^{-s\phi(x,t)}\left[(q_{1} - q_{2})v\right](x,t) - is\left(\partial_{t}\phi(x,t) \right)e^{-s\phi(x,t)}r(x,t) 
+ s^{2}\left|\nabla \phi(x,t)\right|^{2}e^{-s\phi(x,t)}r(x,t) \\ -\; 2se^{-s\phi(x,t)}\left(\nabla \phi \cdot \nabla r \right)(x,t) - s\Delta \phi(x,t) e^{-s\phi(x,t)}r(x,t),
\end{split}
\end{equation}
we note from \eqref{Partial_t_Second_Order_Difference} that
\begin{equation} \label{weighted_Partial_t_Second_Order_Difference}
\begin{cases}
\left(i\partial_t + \Delta + p(x)\right)r_{1}(x,t) = m(x,t) \text{ on } \Omega \times (-T, T),\\
r_{1}(x,t) = 0 \text{ on } \partial \Omega \times (-T, T),\\
r_{1}(x,0) = ie^{-s\phi(x,0)}\sum_{m = 0}^{k}\binom{k}{m} \cdot \left[q_1(x) - q_2(x)\right] \cdot
\partial^{m}_{z_1}\partial^{k - m}_{z_2}N(0,0) \cdot \left[f(x)\right]^{k} \text{ on } \Omega \times \{0\}.
\end{cases}
\end{equation}

\indent The following estimate is proved in Lemma $3.1$ of \cite{MainMethod}:
\begin{equation} \label{initial_data_estimate_0_0}
\Vert r_{1}(\cdot\,,0) \Vert^{2}_{L^{2}(\Omega_{1})} \leq 
C\left(\Vert r_{1}\Vert^{2}_{L^{2}(-2T/3, 2T/3; H^{1}(\Omega_{2}))} + \int_{-2T/3}^{2T/3}\int_{\Omega_{2}}|m(x,t)r_{1}(x,t)|\,dx\,dt \right)
\end{equation}
for some constant $C = C(T) > 0$. Further, by applying the Cauchy-Schwarz Inequality upon expressing ``$e^{-2s\phi(x,t)} = \left(s e^{-s\phi(x,t)}\right)\left(s^{-1} e^{-s\phi(x,t)}\right)$" and ``$se^{-2s\phi(x,t)} =\left(e^{-s\phi(x,t)}\right)\left(s e^{-s\phi(x,t)}\right) $" wherever convenient, we obtain
\begin{equation} \label{initial_data_estimate_0}
\begin{split}
\Vert e^{-s\phi(\cdot\,,0)}r(\cdot\,,0) \Vert^{2}_{L^{2}(\Omega_{1})} \leq 
C\bigg(\frac{1}{s^{2}}\Vert e^{-s\phi}(q_{1} - q_{2}) \Vert^{2}_{L^{2}(\Omega_{2} \times (-2T/3, 2T/3))}\Vert v \Vert_{L^{\infty}(-2T/3, 2T/3; H^{2}(\Omega_{2}))}\\ +\; s^{2}\Vert e^{-s\phi}r \Vert^{2}_{L^{2}(\Omega_{2} \times (-2T/3, 2T/3))}
+\;  \Vert e^{-s\phi}\nabla r \Vert^{2}_{L^{2}(\Omega_{2} \times (-2T/3, 2T/3))}\bigg),
\end{split}
\end{equation}
where $ s \geq s_{0}$ for a constant $s_{0} = s_{0}(\Omega_{2}) > 1$ and $C = C(\Omega_{2}, T) > 0$. We refer the reader to Lemma $3.2$ of \cite{MainMethod} for a similar computation. Using \eqref{weighted_Partial_t_Second_Order_Difference}, Lemma \ref{L2_H2_Partial_t_Second_Order_Difference_estimates} to estimate $\Vert v \Vert_{L^{\infty}(-2T/3, 2T/3; H^{2}(\Omega_{2}))}$, and the fact that $e^{-s\phi(x,t)} \leq e^{-s\phi(x,0)}$, we further obtain
\begin{equation} \label{initial_data_estimate_1_0}
\begin{split}
\Vert e^{-s\phi(\cdot\,,0)}(q_{1} - q_{2})f^{k} \Vert^{2}_{L^{2}(\Omega_{1})} \leq 
C_{3}\bigg(\frac{1}{s^{2}}\Vert e^{-s\phi(\cdot, 0)}(q_{1} - q_{2})\Vert^{2}_{L^{2}(\Omega_{2} \times (-2T/3, 2T/3))}\\ +\; s^{2}\Vert e^{-s\phi}r \Vert^{2}_{L^{2}(\Omega_{2} \times (-2T/3, 2T/3))}
+\;  \Vert e^{-s\phi}\nabla r \Vert^{2}_{L^{2}(\Omega_{2} \times (-2T/3, 2T/3))}\bigg)
\end{split}
\end{equation}
for some constant $C_{3} = C_{3}(\Omega_{2}, T, \gamma_{+}) > 0$. Taking $s = s(C_{3}) \in [s_{0}, \infty)$ sufficiently large to absorb the first RHS term into LHS while using the assumption that $|f(x)| \geq \gamma_{-}$ for $x \in \Omega \setminus \omega$ yields
\begin{equation} \label{initial_data_estimate_1_0_0}
\Vert e^{-s\phi(\cdot\,,0)}(q_{1} - q_{2})\Vert^{2}_{L^{2}(\Omega \setminus \omega)} \leq 
C\left(s^{2}\Vert e^{-s\phi}r \Vert^{2}_{L^{2}(\Omega_{2} \times (-2T/3, 2T/3))}
+ \Vert e^{-s\phi}\nabla r \Vert^{2}_{L^{2}(\Omega_{2} \times (-2T/3, 2T/3))}\right).
\end{equation}
Using the estimate $\phi(x,0)  \leq e^{2\lambda\Vert \psi\Vert_{\infty}}/T_{1}^{2}$ in a manner similar to that of \eqref{initial2_stability_est_0}-\eqref{initial2_stability_est} in Section \ref{main_theorem_proof_geometric_condition}, we obtain
\begin{equation} \label{initial_data_estimate_1}
\Vert  q_{1} - q_{2} \Vert^{2}_{L^{2}(\Omega \setminus \omega)} \leq C\left(s^{2}\Vert e^{-s\phi}r \Vert^{2}_{L^{2}(\Omega_{2} \times (-2T/3, 2T/3))}
+ \Vert e^{-s\phi}\nabla r \Vert^{2}_{L^{2}(\Omega_{2} \times (-2T/3, 2T/3))} \right)
\end{equation}
for some constants $s \geq s_{0}$ and $C = (\Omega,\omega, T, M, \gamma_{-}, \gamma_{+}, \lambda) > 0$. We have omitted the dependence of $C$ on $s = s(C_{3})$ since $C_{3}$ itself depends on the parameters $\Omega_{2}, T, \gamma_{+}$. \\ \indent
Now, after applying the Carleman estimate in Proposition \ref{standard_Schrodinger_Carleman_estimate_0} (with $T_{1} = 2T/3$, $\lambda \geq \text{max}\{\lambda_{*},1\}, s \geq \max\{s_{*}, s_{0}\}$), the estimates in Lemma \ref{L2_H2_Partial_t_Second_Order_Difference_estimates}, and the unique continuation estimate \eqref{logarithmic_observability_inequality_statement}, we obtain
\begin{equation}\label{initial_data_estimate_2}
\begin{split}
\Vert  q_{1} - q_{2} \Vert^{2}_{L^{2}(\Omega \setminus \omega)} \leq \frac{C}{s}\Vert  q_{1} - q_{2} \Vert^{2}_{L^{2}(\Omega \setminus \omega)} + Cs^{2}\left(\frac{1}{\gamma} + e^{-\mu \gamma} +  e^{\mu \gamma}\int_{\Gamma}\int_{-T}^{T}  \left|\frac{\partial r}{\partial \nu}(x,t)\right|^{2} \,dt \,dS\right)
\end{split}
\end{equation}
for $s \geq \max\{s_{*}, s_{0}\}$, $\gamma \geq \gamma_{*}$, and constants $C = C(\Omega, \omega, M, \tilde{M}, \gamma_{-}, \gamma_{+}, T, T_{0}) > 0$, $\mu = \mu(\omega, T, T_{0}) > 0$. \\ \indent
Recall in relation to \eqref{stability_estimate_step_3_0_0} that we have $\gamma_{*} = \sigma_{*}(\omega)/T_{0}$. So, there exists a constant $\delta_{*} = \delta_{*}(\omega, T, T_{0}) > 0$ such that
\begin{equation}
\left\Vert \frac{\partial r}{\partial \nu}\right\Vert_{L^{2}(\Gamma \times (-T, T))} < \delta_{*}
\end{equation}
implies
\begin{equation}
\gamma = \frac{1}{\mu}\ln \left\Vert \frac{\partial r}{\partial \nu}\right\Vert_{L^{2}(\Gamma \times (-T, T))}^{-1} \geq \gamma_{*}.
\end{equation}
Taking $s \geq \max\{s_{*}, s_{0}\}$ sufficiently large in \eqref{initial_data_estimate_2}, we obtain
\begin{equation} \label{partial_data_stability_estimate}
\Vert  q_{1} - q_{2} \Vert_{L^{2}(\Omega \setminus \omega)} \leq C\left[\left(\ln \left\Vert \frac{\partial r}{\partial \nu} \right\Vert^{-1}_{L^{2}(\Gamma \times (-T, T))} \right)^{-1} + \left\Vert \frac{\partial r}{\partial \nu} \right\Vert_{L^{2}(\Gamma \times (-T, T))} \right]^{1/2}
\end{equation}
for a constant $C = C(\Omega, \omega, M, \tilde{M}, \gamma_{-}, \gamma_{+}, T) > 0$ where we have omitted ``$T_{0}$" due to its dependence on $\omega$.\\ \indent
On the other hand, if we have that 
\begin{equation}
\left\Vert \frac{\partial r}{\partial \nu} \right\Vert_{L^{2}(\Gamma \times (-T, T))} \geq \delta_{*},
\end{equation}
Lemma \ref{L2_H2_Partial_t_Second_Order_Difference_estimates} implies that we can write
\begin{equation} \label{lower_bound_Neumann_data_stability_estimate_step_5}
\int_{\Omega_{2}\setminus \Omega_{1}}\int_{-2T/3}^{2T/3}\left|\nabla r(x,t) \right|^{2} + \left|r(x,t) \right|^{2}\,dt\,dx \leq \frac{C}{\delta_{*}}\delta_{*} 
\leq C\left\Vert \frac{\partial r}{\partial \nu} \right\Vert_{L^{2}(\Gamma \times (-T, T))}
\end{equation}
for some $C = C(\Omega, M, T) > 0$ instead of \eqref{logarithmic_observability_inequality_statement}. Starting from \eqref{initial_data_estimate_1} and applying the Carleman estimate in Proposition \ref{standard_Schrodinger_Carleman_estimate_0}, we obtain
\begin{equation} \label{partial_data_stability_estimate_1}
\Vert  q_{1} - q_{2} \Vert_{L^{2}(\Omega \setminus \omega)} \leq C\left\Vert \frac{\partial r}{\partial \nu} \right\Vert^{1/2}_{L^{2}(\Gamma \times (-T, T))}
\end{equation}
for some constant $C = C(\Omega, M, \tilde{M}, \gamma_{-}, \gamma_{+}, \delta_{*},T) > 0$. \\ \indent
We pass from the interval $(-T, T)$ to $(0,T)$ using Remark \ref{symmetric_extension_remark} and this completes the proof of \eqref{main_theorem_partial_stable_determination} upon recalling from \eqref{derivative_interchange} that one can interchange the order of differentiation in $\partial r/\partial \nu$. We similarly prove \eqref{main_theorem_partial_stable_determination_0} using $r_{1}(x,t)$ defined by \eqref{Partial_t_First_Linear_Difference} and estimates similar to those in Lemma \eqref{L2_H2_Partial_t_Second_Order_Difference_estimates}. Combining Remark \ref{derivative_interchange_and_main_theorem_statement} with \eqref{main_theorem_partial_stable_determination_0}-\eqref{main_theorem_partial_stable_determination} yields unique determination of $p_{j}(x),\,q_{j}(x)\;(j = 1, 2)$ assuming knowledge of the coefficients near $\partial \Omega$.

\appendix
\section{The Solutions to \eqref{FirstLinear_for_p1} and \eqref{Partial_t_Second_Order_Difference} on $(-T, T)$} \label{Extension_Discussion}
We wish to check that the solutions to equations \eqref{FirstLinear_for_p1} and \eqref{Partial_t_Second_Order_Difference} are such that ${^{(1)}}u_{\epsilon}(x,t) = \overline{{^{(1)}}u_{\epsilon}(x,-t)}$ and $r(x,t) = -\overline{r(x,-t)}$ when $t \in (-T,0)$. For example, let us verify that for $t \in (-T,0)$, $r(x,t) = -\overline{r(x,\tau)}$ with $\tau = -t \in (0,T)$ satisfies \eqref{Partial_t_Second_Order_Difference}. That is, we seek to show that
\begin{align*}
\left(i\partial_t + \Delta + p(x)\right) r(x,t) = -\sum_{m = 0}^{k}\binom{k}{m} \cdot \left[q_1(x) - q_2(x)\right] \cdot 
\partial^{m}_{z_1}\partial^{k - m}_{z_2}N(0,0) \cdot\\ \bigg[m \cdot \partial_{t} \left({^{(1)}}u_{\epsilon}(x,t)\right) \cdot \left({^{(1)}}u_{\epsilon}(x,t)\right)^{m-1}  \cdot  \Big(\overline{^{(1)}u_{\epsilon}(x,t)}\Big)^{k - m} + \\
(k - m) \cdot \partial_{t}\left(\overline{{^{(1)}}u_{\epsilon}(x,t)}\right) \cdot  \left({^{(1)}}u_{\epsilon}(x,t)\right)^{m} \cdot 
\left(\overline{{^{(1)}}u_{\epsilon}(x,t)}\right)^{k - m - 1}
\bigg]
\end{align*}
where ${^{(1)}}u_{\epsilon}(x,t) = \overline{{^{(1)}}u_{\epsilon}(x,\tau)}$ for $t \in (-T,0)$. We may directly compute
\begin{align*}
i\partial_t r(x,t) + \Delta r(x,t) + p(x)r(x,t) \\  =
i\partial_t (-\overline{r(x,\tau)}) + \Delta (-\overline{r(x,\tau)})+ p(x)(-\overline{r(x,\tau)}) \\ =
i\partial_{\tau} \overline{r(x,\tau)} - \Delta \overline{r(x,\tau)} - p(x)\overline{r(x,\tau)}.
\end{align*}
Now, suppose that $q_1(x)-q_{2}(x)$ and
\begin{equation*}
\partial^{m}_{z_{1}}\partial^{k - m}_{z_{2}}N(0,0)\;(m = 0,\ldots,k)
\end{equation*}
are real-valued. Then, we have from conjugating and negating \eqref{Partial_t_Second_Order_Difference} that
\begin{align*}
\left(i\partial_{\tau} - \Delta  - p(x)\right)\overline{r(x,\tau)} = \sum_{m = 0}^{k}\binom{k}{m} \cdot \left[q_1(x) - q_2(x)\right] \cdot 
\partial^{m}_{z_1}\partial^{k - m}_{z_2}N(0,0) \cdot\\ \bigg[m \cdot \partial_{\tau} \left(\overline{{^{(1)}}u_{\epsilon}(x,\tau)}\right) \cdot \left(\overline{{^{(1)}}u_{\epsilon}(x,\tau)}\right)^{m-1}  \cdot  \Big({^{(1)}}u_{\epsilon}(x,\tau)\Big)^{k - m} + \\
(k - m) \cdot \partial_{\tau} \left({^{(1)}}u_{\epsilon}(x,
\tau)\right) \cdot  \left(\overline{{^{(1)}}u_{\epsilon}(x,\tau)}\right)^{m} \cdot 
\left({^{(1)}}u_{\epsilon}(x,\tau)\right)^{k - m - 1}
\bigg].
\end{align*}
The required verification follows from the fact that ${^{(1)}}u_{\epsilon}(x,t) = \overline{{^{(1)}}u_{\epsilon}(x,\tau)}$ and that $\partial_{\tau} \overline{^{(1)}u_{\epsilon}(x,\tau)} = -\partial_{t} \overline{^{(1)}u_{\epsilon}(x,\tau)} = -\partial_{t} {^{(1)}}u_{\epsilon}(x,t)$. One similarly verifies that if $f(x)$ is real-valued, we can write ${^{(1)}}u_{\epsilon}(x,t) = \overline{{^{(1)}}u_{\epsilon}(x,-t)}$.

\end{document}